\documentclass[11pt]{amsart}
\usepackage{amsmath,amsthm, amscd, amssymb, amsfonts,mathtools}
\usepackage[all]{xy}
\usepackage{mathrsfs}
\usepackage[inline]{enumitem}
\usepackage{hyperref}
\usepackage{multicol}
\usepackage{rotating}

\usepackage[ansinew]{inputenc}
\usepackage{graphicx,fancyhdr}

\usepackage[dvips, dvipsnames, usenames]{color}

\newcommand{\xrightarrowdbl}[2][]{%
\xrightarrow[#1]{#2}\mathrel{\mkern-14mu}\rightarrow
}

\newcommand{\yd}[1]{{}^{ #1 }_{ #1 }\mathcal{YD}}

\newcommand{\comment}[1]{}

\newcommand{\fd}{finite-dimensional}
\newcommand{\trid}{\triangleright}
\newcommand{\Gx}{\mathbb G_X}

\numberwithin{equation}{section}
\newtheorem{theorem}{Theorem}[section]
\newtheorem{lemma}[theorem]{Lemma}
\newtheorem{coro}[theorem]{Corollary}

\newtheorem{prop}[theorem]{Proposition}

\theoremstyle{definition}
\newtheorem{definition}[theorem]{Definition}
\newtheorem{definitionlemma}[theorem]{Definition--Lemma}
\newtheorem{example}[theorem]{Example}
\newtheorem{question}{Question}

\newtheorem{assumption}[theorem]{Assumption}

\theoremstyle{remark}
\newtheorem{remark}[theorem]{Remark}

\newcommand\pf{\begin{proof}}
\newcommand\epf{\end{proof}}

\newcommand{\ku}{ \Bbbk}
\newcommand{\kut}{ \ku^{\times}}

\newcommand{\G}{\mathbb G}
\newcommand{\gb}{\mathbf g}
\newcommand{\chib}{\boldsymbol \chi}

\newcommand{\ct}{\mathtt{c}}
\newcommand{\FK}{\mathtt{FK}}

\newcommand{\I}{\mathbb I}

\newcommand\N{\mathbb N}

\newcommand{\bq}{\mathfrak{q}}
\newcommand\Sb{\mathbb S}
\newcommand{\sn}{\Sb_n}

\newcommand\Z{\mathbb Z}

\renewcommand{\_}[1]{_{\left( #1 \right)}}
\renewcommand{\^}[1]{^{\left( #1 \right)}}

\newcommand\cA{\mathcal{A}}
\newcommand\cB{\mathcal{B}}

\newcommand\dpn{\widetilde{\mathcal{B}}}

\newcommand{\Cr}{\mathscr{C}}

\newcommand{\E}{\mathcal{E}}
\newcommand{\F}{\mathcal F}
\newcommand{\cH}{\mathcal{H}}
\newcommand{\hcenter}{\mathcal{HZ}}
\newcommand{\hcocenter}{\mathcal{HC}}

\newcommand{\cJ}{\mathcal{J}}

\newcommand{\Oc}{{\mathcal O}}

\newcommand{\Ss}{{\mathcal S}}

\newcommand{\g}{\mathfrak g}

\newcommand{\ug}{\mathfrak u}
\newcommand{\ag}{\mathfrak a}

\newcommand{\Mg}{\mathfrak M}
\newcommand{\Ig}{\mathfrak I}

\newcommand{\alg}{\cH}

\newcommand{\ad}{\operatorname{ad}}
\newcommand{\Alg}{\operatorname{Alg}}

\newcommand{\End}{\operatorname{End}}

\newcommand{\Imm}{\operatorname{Im}}

\newcommand{\id}{\operatorname{id}}
\newcommand{\gr}{\operatorname{gr}}

\newcommand{\pd}{\operatorname{pd}}
\newcommand{\gldim}{\operatorname{gldim}}
\newcommand{\lgldim}{\operatorname{l.gldim}}
\newcommand{\rgldim}{\operatorname{r.gldim}}
\newcommand{\corad}{\operatorname{corad}}

\newcommand{\GK}{\operatorname{GKdim}}
\newcommand{\Hom}{\operatorname{Hom}}
\newcommand\ord{\operatorname{ord}}
\newcommand{\rep}{\operatorname{rep}}

\newcommand{\Res}{\operatorname{Res}}
\newcommand{\Reg}{\operatorname{Reg}}
\newcommand{\Specmax}{\operatorname{Specmax}}
\newcommand{\supp}{\operatorname{supp}}

\newcommand{\toba}{\mathscr{B}}
\newcommand{\ot}{\otimes}

\newcounter{tabla}\stepcounter{tabla}

\begin{document}

\title[Finite-by-cocommutative  Hopf algebras]{A class of finite-by-cocommutative  Hopf algebras}

\author[Andruskiewitsch, Natale, Torrecillas]
{Nicol\'as Andruskiewitsch,  Sonia Natale, Blas Torrecillas}

\address{ N.~A., S.~N.: FaMAF-CIEM (CONICET), Universidad Nacional de C\'ordoba,
Me\-dina A\-llen\-de s/n, Ciudad Universitaria, 5000 C\' ordoba,  Argentina.} 

\email{nicolas.andruskiewitsch@unc.edu.ar, natale@famaf.unc.edu.ar}

\address{B. T.: Universidad de Almer\'\i a, Dpto. \'Algebra y An\'alisis Matem\'atico.
E04120 Almer\'\i a, Spain}
\email{btorreci@ual.es}

\thanks{\noindent 2020 \emph{Mathematics Subject Classification.}
16T05; 18M05. \newline 
 The research on the subject of this paper was started during a visit of N. A.  and S. N.
 to the University of Almer\'\i a, whose hospitality is gratefully acknowledged.
 The work of N. A. was  partially supported by CONICET (PIP 11220200102916CO),
FONCyT-ANPCyT (PICT-2019-03660) and Secyt (UNC). The research of S. N. was partially supported by Secyt (UNC).
The research of B. T.  was supported partially by PID2020-113552GB-I00 from MICIN, FEDER- UAL18-FQM-B042-A and P20-0770 from Junta de Andalucía}

\begin{abstract}
We present a rich source of Hopf algebras starting from a cofinite central extension of a Noetherian Hopf algebra
and a subgroup of the algebraic group of characters of the central Hopf subalgebra. The construction is transparent
from a Tannakian perspective. We determine when the new 
Hopf algebras are co-Frobenius, or cosemisimple, or Noetherian, or regular, or have finite Gelfand-Kirillov dimension.
\end{abstract}

\maketitle

\setcounter{tocdepth}{1}

\section{Introduction}\label{section:introduction}

A  Hopf algebra $H$ is \emph{commutative-by-finite} if  it has a normal Hopf subalgebra $A$ such that 
 $A$ is commutative and $H$ is a finitely generated $A$-module.
In other words there is an exact sequence   of Hopf algebras 
$\xymatrix{ A  \ar@{^{(}->}[r]  & H \ar@{->>}[r]  & \ug}$
where $A$ is commutative and $\ug$ is finite-dimensional.
 There are various remarkable families of commutative-by-finite Hopf algebras arising from the theory of quantum groups.
A systematic study of \emph{affine} commutative-by-finite Hopf algebras was started in \cite{bc}; here  `affine' means that $H$ is a finitely generated algebra.

\medbreak
A Hopf algebra $K$ is \emph{finite-by-cocommutative} if  it fits into an exact sequence of Hopf algebras  
$\xymatrix{ \ag  \ar@{^{(}->}[r]  & K \ar@{->>}[r]  & U}$ where $U$ is cocommutative and $\ag$ is finite-dimensional.
Lusztig's quantum groups at roots of 1 are finite-by-cocommutative.
An example of a finite-by-cocommutative Hopf algebra appeared in \cite{ace} to disprove a conjecture on co-Frobenius Hopf algebras. 
A family of examples
containing that one and characterized by suitable properties was presented in \cite{li-liu}. The finite dual $H^{\circ}$
of an affine commutative-by-finite Hopf algebra $H$ was studied in \cite{bcj}; it is finite-by-cocommutative. 

\medbreak
The goal of  this paper is to present and study a family  of  finite-by-cocommutative Hopf algebras.
We assume that the base field $\ku$ is algebraically closed and has characteristic 0.
To start with, consider a Noetherian Hopf algebra $H$ with a central Hopf subalgebra $A$;
 we set $H_{\varepsilon} = H/  HA^{+}$, where $A^{+} = \ker \varepsilon_{\vert A}$. 
Thus  we have  an extension of Hopf algebras   
\begin{align}
\tag{$\E$}  
A  \hookrightarrow H  \twoheadrightarrow  H_{\varepsilon}.
\end{align}

In Section \ref{sec:hopf-systems}, for any subgroup $\varGamma$ of the 
pro-affine algebraic group $G = \Alg(A, \ku)$, we define a suitable subgroup $\varGamma_{\hspace{-2pt}\mathrm{fd}}$ of $\varGamma$ and a Hopf subalgebra $\alg(\varGamma)$ of the finite dual 
$H^{\circ}$, which  is an  extension of Hopf algebras
\begin{align}
\tag{$\F^{\varGamma}$} H_{\varepsilon}^{\circ}   \hookrightarrow  \alg(\varGamma) \twoheadrightarrow \ku \varGamma_{\hspace{-2pt}\mathrm{fd}}.
\end{align} 
From a Tannakian perspective, 
the category of finite-dimensional comodules over $\alg(\varGamma)$  is equivalent to
 the full subcategory $\Cr_{\varGamma}$  of  the category $\rep H$ of finite-dimensional $H$-modules
 such that  the action of $A$ is semisimple and by characters in $\varGamma$; thus
 the objects of $\Cr_{\varGamma}$ bear a $\varGamma$-grading.

\medbreak
In Section \ref{sec:finite-by-coco} we assume further that the extension $(\E)$ is cleft and that 
$\dim H_{\varepsilon}< \infty$. Then $H$ is a finitely generated $A$-module, $A$ is Noetherian, 
$G$ is an algebraic group, $\varGamma_{\hspace{-2pt}\mathrm{fd}} = \varGamma$,   $(\F^{\varGamma})$ becomes
$H_{\varepsilon}^{*} \hookrightarrow  \alg(\varGamma)  \twoheadrightarrow \ku \varGamma$ and so
 $\alg(\varGamma)$ is finite-by-cocommutative.
We establish several properties of $\alg(\varGamma)$:

\begin{theorem}\label{thm:main}
\begin{enumerate}[leftmargin=*,label=\rm{(\roman*)}] 
\item If $\varGamma$ is finitely-generated, then $\alg(\varGamma)$ is affine. 

\medbreak
\item $\alg(\varGamma)$ is co-Frobenius.

\medbreak
\item $\alg(\varGamma)$ is cosemisimple if and only if $H_{\varepsilon}$ is semisimple.

\medbreak
\item If $\varGamma$ is finitely-generated, then
$\GK\alg(\varGamma) < \infty$   if and only if  $\varGamma$  is nilpotent-by-finite.

\medbreak
\item $\alg(\varGamma)$ is Noetherian if and only if $\varGamma$ is polycyclic-by-finite.

\medbreak
\item If  $\alg(\varGamma)$ is Noetherian, then it is regular iff $H_{\varepsilon}$ is semisimple.
\end{enumerate}
\end{theorem}
Theorem \ref{thm:main} is proved  in in Section \ref{sec:finite-by-coco}, see
Remark \ref{rem:H(Gamma)affine} and Theorems \ref{thm:cofrob}, \ref{thm:coss},
\ref{thm:gkdim}, \ref{thm:Noetherian} and \ref{thm:regular}.
The keys to these results are that $\alg(\varGamma)$ is \emph{strongly} $\varGamma$-graded
and that $\varGamma$, being a subgroup of an algebraic group, is linear, i.e.,  embedable into $GL(n, \ku)$ for some $n$.

Along the way to Theorem \ref{thm:regular} we found a general, apparently new, criterium: 
A Hopf algebra that contains a finite-dimensional normal non-semisimple Hopf subalgebra  is not regular, see Lemma \ref{lemma:regular-gral}.

The contents of the paper are organized in the following way. 
Section \ref{section:Preliminaries} contains expositions of known facts needed along the article.
In Section \ref{sec:hopf-systems} we present the Hopf algebras $\alg(\varGamma)$ and the extensions
$(\F^{\varGamma})$ in a general context and establish some basic properties.
In Section \ref{sec:finite-by-coco} we study the Hopf algebras $\alg(\varGamma)$ under the restrictions above.
Section \ref{section:examples} is devoted to  examples.

\subsection*{Notations}
The natural  numbers are denoted by $\N$, and $\N_0=\N\cup \{0\}$. 
Given $m< n \in\N_0$, we set $\I_{m,n}=\{i \in \N_0: m \leq i \leq m\}$ and $\I_{n}=\I_{1, n}$.
`Algebra' means associative unital algebra. The space of algebra homomorphism from a $\ku$-algebra $A$ to a  $\ku$-algebra $B$ is denoted by $\Alg(A, B)$.
The category of  finite-dimensional left $R$-modules, where $R$ is an  algebra, is denoted by $\rep R$.
All Hopf algebras  are supposed to have bijective antipode. We write $M\leq N$ to express that 
$M$ is a subobject of $N$ in a given category. The notation for Hopf algebras is standard: $\Delta$ is the comultiplication,
$\varepsilon$ is the counit, $\Ss$ is the antipode. For the comultiplication and the coactions we use the Heynemann-Sweedler notation.

\section{Preliminaries}\label{section:Preliminaries}
We refer the reader to \cite{montgomery,rad-libro} for the basic facts about Hopf algebras used throughout the paper.
Given a coalgebra $K$ and  an algebra $H$, the
group of invertible elements in $\Hom(K, H)$ with respect to the convolution  is denoted by $\Reg(K,H)$.

\subsection{Cleft comodule algebras}
These  were studied in \cite{bcm,doi-takeuchi}; we recall the relevant facts.
Let $H$ be a Hopf algebra, let $R$ be a right comodule algebra with coaction $\rho: R \to R \otimes H$
and let $R^{\operatorname{co} H} =\{x\in R: \rho(x) = x \otimes 1\}$. 
One defines similarly ${}^{\operatorname{co} H}T$ for a left comodule algebra $T$.

\medbreak
For instance, if $\pi: C\to B$ is a  Hopf algebra map, then $C$ is a right, respectively left, comodule algebra
via $\rho = (\id \otimes \pi)\Delta$, resp.  $\lambda = (\pi \otimes \id)\Delta$. 
The  algebras of right and left coinvariants of $\pi$ are
\begin{align*}
C^{\operatorname{co} \pi} & = C^{\operatorname{co} B} =\{x\in C: (\id \otimes \pi)\Delta(x) = x \otimes 1\}, 
\\
{}^{\operatorname{co} \pi} C &=  {}^{\operatorname{co} B}C = \{x\in C: (\pi \otimes \id)\Delta(x) = 1 \otimes x\}.
\end{align*}

\medbreak
We consider three  properties of the extension of algebras $R^{\operatorname{co} H} \subset R$:

\medbreak
\begin{enumerate}[leftmargin=*,label=\rm{(\roman*)}] 
\item\label{item:galois}   $R^{\operatorname{co} H} \subset R$ is  $H$-Galois, if the canonical map
$\operatorname{can}:R\otimes_{R^{\operatorname{co} \pi}} R \to R \otimes H$, given by
$x\otimes y \mapsto (x\otimes 1)\rho(y)$, is bijective.

\medbreak
\item\label{item:normalbasis} 
$R^{\operatorname{co} H} \subset R$
has the  normal basis property if $R \simeq R^{\operatorname{co} H} \otimes H$ as left $R^{\operatorname{co} H}$-modules and
right $H$-comodules.

\medbreak
\item\label{item:cleft}  
$R^{\operatorname{co} H} \subset R$ is cleft if there exists $\chi\in \Reg( H, R)$   such that $\chi$ is a morphism of $H$-comodules, 
i.e., $\rho\chi = (\chi\otimes \id)\Delta$.
\end{enumerate}

\begin{theorem}\label{thm:doi-takeuchi}\cite{doi-takeuchi} 
The extension $R^{\operatorname{co} H} \subset R$ is cleft if and only if it is $H$-Galois and has  the  normal basis property. \qed
\end{theorem}

\begin{example}\label{exa:strongly-graded}  \cite{ulbrich}
Let $G$ be a group with unit $e$ and let $R$ be an algebra.  A  $\ku G$-comodule algebra structure on $R$
is the same as a $G$-grading  of algebras $R = \oplus_{g \in G} R_g$; here $R^{\operatorname{co} H} = R_{e}$. Now 
`$R$ is strongly $G$-graded' means that 
\begin{align*}
R_gR_h &=  R_{gh},& \text{ for all  }g, h &\in G. 
\end{align*}
Then $R_{e} \subset R$ is $\ku G$-Galois if and only if $R$ is strongly $G$-graded.
\end{example}

\subsection{Extensions of Hopf algebras}\label{subsec:extensions}
This notion was considered in many papers. see e.g. \cite{kac-extension,ad, sch,hofstetter,singer,takeuchi,majid,byott}.
Following \cite{a-can,ad} together with \cite{sch} we say that the sequence of morphisms of Hopf algebras
\begin{align*}
A\xhookrightarrow[]{\iota} C \xrightarrowdbl[]{\pi} B
\end{align*}
is exact if the following conditions hold:
\begin{multicols}{2}
\begin{enumerate}[leftmargin=*,label=\rm{(\roman*)}]
\item\label{suc-exacta-1} $\iota$ is injective.
\item\label{suc-exacta-2} $\pi$ is surjective.
\item\label{suc-exacta-3} $\ker\pi = C\iota(A)^+$.
\item\label{suc-exacta-4} $\iota(A) = C^{\operatorname{co} \pi}$.
\end{enumerate}
\end{multicols}
In this case we also say that  $C$ is an extension of $B$ by $A$.

\medbreak 
The left and right adjoint actions of $C$, denoted by $\ad_{\ell}, \ad_{r}: C \to \End C$, are given by
\begin{align*}
\ad_{\ell}(x)(y) &= x\_{1}y \Ss (x\_{2}), &
\ad_{r}(x) (y)&= \Ss (x\_{1}) y x\_{2}, & x, y&\in C.
\end{align*}
Notice that
\begin{align*}
\ad_{r}(x) (y)&= \Ss^{-1} \left(\ad_{\ell}(\Ss (x))\Ss (y)\right), & x, y&\in C.
\end{align*}
A Hopf subalgebra $A$ of a Hopf algebra $C$ is \emph{normal} if it stable under one of, hence  both,   the adjoint actions.

\begin{lemma}\label{lemma:exact-sequence-hopf} Let  $A$ be Hopf subalgebra of a Hopf algebra $C$and $B\coloneqq C/A^+C$.

\begin{enumerate}[leftmargin=*,label=\rm{(\roman*)}] 
\item\label{item:exact-sequence-hopf-normal} If $A$ is normal, then $A^+C = CA^+$. 

\medbreak
\item\label{item:exact-sequence-hopf-normal2}
 If $A^+C = CA^+$, then it  is a Hopf ideal,  $B$ is a Hopf algebra and the quotient map $\pi:C \to B$
is a morphism of Hopf algebras. 

\medbreak
\item\label{item:exact-sequence-hopf-fflat}  If $A^+C = CA^+$ and $C$ is a faithfully flat $A$-module (under left or right multiplication),
then $A\xhookrightarrow[]{\iota} C \xrightarrowdbl[]{\pi} B$ is exact, $\pi$ is faithfully coflat
 and $A$ is normal.
\end{enumerate}
\end{lemma}

The converse in \ref{item:exact-sequence-hopf-normal} and whether Hopf algebras are faithfully flat over Hopf subalgebras
are open questions.

\pf \ref{item:exact-sequence-hopf-normal}, \ref{item:exact-sequence-hopf-normal2} are easy; see \cite[3.4.3]{montgomery} for
\ref{item:exact-sequence-hopf-fflat}; cf. \cite[1.2.4]{ad}, 
\cite[1.4]{sch} \cite{takeuchi-quotient}.
\epf

\begin{remark}\label{rem:exact-sequence-hopf-fflat}
A Hopf algebra $C$ is faithfully flat over a Hopf subalgebra  $A$  provided that either of the following  conditions hold:

\begin{enumerate}[leftmargin=*,label=\rm{(\roman*)}] 

\item\label{item:exact-sequence-hopf-finite}  $\dim C < \infty$, in fact $C$ is a free  $A$-module \cite{nz};

\medbreak
\item\label{item:exact-sequence-hopf-fflat2}\cite[3.3]{sch}   $A$ is  central and $C$ is Noetherian;

\medbreak
\item\cite[2.1]{sch}  $A$ is normal and $\dim A < \infty$; in fact $C$ is a free  $A$-module. 
\end{enumerate}
\end{remark}

\medbreak Let $C$ be a Hopf algebra.  The left and right coadjoint actions of $C$ are the (left and right) comodule structures
$\varrho_{\ell}, \varrho_{r}: C \to C \otimes C$ given by
\begin{align*}
\varrho_{\ell}(x) &= x\_{1} \Ss (x\_{3}) \otimes x\_{2}, &
\varrho_{r}(x) &= x\_{2} \otimes \Ss (x\_{1})  x\_{3}, & x&\in C.
\end{align*}
Let $\tau$ be the usual flip. Notice that
\begin{align*}
\varrho_{r}(x) &= \left(\Ss^{-1} \otimes \Ss^{-1}\right) \tau \left(\varrho_{\ell}(\Ss (x))\right), & x&\in C.
\end{align*}

A surjective Hopf algebra map $\pi: C \to B$ is \emph{conormal}, or simply $B$ is a conormal quotient of $C$,
if $\ker \pi$ is a subcomodule for one of, hence both,  $\varrho_{\ell}$ and $\varrho_{r}$ (notice a change of terminology with respect to \cite{ad}).

\begin{lemma}\label{lemma:exact-sequence-hopf2} Let  $\pi: C \to B$ be a surjective morphism of Hopf algebras.

\begin{enumerate}[leftmargin=*,label=\rm{(\roman*)}] 
\item\label{item:exact-sequence-hopf-conormal} If $\pi$ is conormal, then $C^{\operatorname{co} \pi} 
= {}^{\operatorname{co} \pi}C$. 

\medbreak
\item\label{item:exact-sequence-hopf-conormal2}
 $A \coloneqq C^{\operatorname{co} \pi}$ equals $ {}^{\operatorname{co} \pi}C$ iff $A$ is a Hopf subalgebra of $C$. 

\medbreak
\item\label{item:exact-sequence-hopf-fcoflat} If $\pi$ is conormal and $C$ is a faithfully coflat  (left or right) $B$-comodule,
then $A\coloneqq C^{\operatorname{co} \pi}\xhookrightarrow[]{\iota} C \xrightarrowdbl[]{\pi} B$ is exact and $\iota$ is
faithfully flat.
\end{enumerate}
\end{lemma}

\pf \ref{item:exact-sequence-hopf-conormal}:  \cite[1.1.7]{ad};
\ref{item:exact-sequence-hopf-conormal2}  \cite[1.1.4]{ad}; 
\ref{item:exact-sequence-hopf-fcoflat}:    \cite[1.2.14]{ad}, \cite[1.4]{sch}, \cite{takeuchi-quotient}.
\epf

\begin{lemma}\label{lema:exact-sequence-hopf-fcoflat}
If  $\pi: C \to B$ is a surjective morphism of Hopf algebras and $B$ is cosemisimple, then 
$C$ is a faithfully coflat $B$-comodule.
\end{lemma}

\pf   $C$ is coflat because $B$ is cosemisimple;
but coflatness for a surjective Hopf algebra map implies faithful coflatness by  \cite[Remark, p. 247]{doi1}.
\epf

\subsection{Cleft extensions of Hopf algebras}\label{subsec:cleft-extensions}
These were considered e.~g. in \cite{a-can,ad,byott,sch2,sch3}.
Given Hopf algebras $H$ and  $K$, we consider the subgroups of $\Reg(K,H)$ given by
\begin{align*}
\Reg_{1} (K,H) &= \{\phi\in \Reg (K,H): \phi(1) = 1\},\\
\Reg_{\varepsilon} (K,H) &= \{\phi \in \Reg (K,H): \varepsilon\phi  = \varepsilon\},\\
\Reg_{1, \varepsilon} (K,H) &= \Reg_{1} (K,H)\cap \Reg_{\varepsilon} (K,H).
\end{align*}

\begin{definitionlemma} \label{def-lemma:extension-hopf}
An exact sequence $A\xhookrightarrow[]{\iota} C \xrightarrowdbl[]{\pi} B$ is \emph{cleft} 
if  satisfies one of the following equivalent conditions:

\begin{enumerate}[leftmargin=*,label=\rm{(\roman*)}]
\item there exists $\chi \in \Reg_{1}(B, C)$ such
that $(\id\otimes \pi)\Delta\chi = (\chi\otimes \id)\Delta$;

\medbreak
\item\label{item:xi}  there exists $\xi \in
\Reg_{\varepsilon}(C, A)$ such that  $ \xi (ac) = a \xi(c), \quad
\forall  a\in A, c \in C$;

\medbreak
\item there exist a morphism of $A$-modules $\xi: C\to A$ and
a  morphism of $B$-comodules $\chi: B\to C$ such that $\xi\chi =
\varepsilon_{B}1_{A}$ and $(\iota\xi) * (\chi\pi) = \id_{C}$.
\end{enumerate}
\end{definitionlemma} 

If this happens, then
 $\xi(1) = 1$, $\varepsilon\chi = \varepsilon$, hence
$\pi\chi = \id_{B}$ and $\xi\iota = \id_{A}$.   
See \cite[3.1.14]{a-can} for a proof.

\subsection{Hopf center and Hopf cocenter} Let $H$ be a Hopf algebra. We recall a few facts from \cite{a-can}.

\begin{itemize}[leftmargin=*] 
\item  \cite[2.2.3]{a-can} There exists a maximal  (with respect to the inclusion)
central Hopf subalgebra $\hcenter(H)$ of $H$, called the  \emph{Hopf center}  of $H$.

\medbreak
\item \cite[2.3.8]{a-can} A quotient Hopf algebra map
$q:H\to K$ is cocentral if it satisfies 
\[(q \otimes \id) \Delta = (q \otimes \id) \Delta^{\operatorname{op}}. \]
There exists a minimal cocentral Hopf algebra map $\bq:  H \to \hcocenter(H)$; by abuse of notation, $\hcocenter(H)$ 
is called the \emph{Hopf cocenter} of $H$.
Here \emph{minimal} means that any cocentral   map $q:H\to K$ factorizes through  $\bq$.

\medbreak
\item If $\dim H < \infty$, then $\hcenter(H)^* \simeq \hcocenter(H^*)$ and $\hcenter(H^*) \simeq \hcocenter(H)^*$.
\end{itemize}  

\begin{lemma}\label{lema:trivial-center}  \cite[3.3.9]{a-can} 
Let  $A\xhookrightarrow[]{\iota} C \xrightarrowdbl[]{\pi} B$ be an exact sequence.
\begin{enumerate}[leftmargin=*,label=\rm{(\roman*)}]
\item If $A$ is central in $C$ and $\hcenter(B) \simeq \ku$, then $A$ is the Hopf center of $C$.

\medbreak 
\item\label{item:trivial-cocenter} If  $\pi: C \to B$ is cocentral and $\hcocenter(A) \simeq \ku$, then $B\simeq \hcocenter(C)$. \qed
\end{enumerate}
\end{lemma}

\begin{remark} Let $H$ be a finite-dimensional simple Hopf algebra (i.e., without proper non-trivial normal Hopf subalgebras).
\begin{enumerate}[leftmargin=*,label=\rm{(\roman*)}]
\item If $H$ is not commutative, then $\hcenter(H) \simeq\ku$.

\medbreak 
\item If $H$ is not cocommutative, then $\hcocenter(H) \simeq\ku$.

\medbreak 
\item (N. A. and H.-J.Schneider, Appendix to \cite{a-can}).
The small quantum groups associated to simple Lie algebras and their parabolic subalgebras are simple Hopf algebras, 
hence their Hopf centers and cocenters are trivial.
\end{enumerate}

\end{remark}

\subsection{Cocycles and twists} \label{subsec:cocycles-twists}
Let $H$ be a Hopf algebra.
A  Hopf 2-cocycle \cite{doi} or simply a cocycle is a convolution invertible inear map $\sigma:H\ot H\to \ku$ that satisfies
$\sigma(x_{(1)},y_{(1)}) \sigma(x_{(2)}y_{(2)},z)  =\sigma(y_{(1)},z_{(1)}) \sigma(x,y_{(2)}z_{(2)})$, 
$\sigma(x,1) = \sigma(1,x)=\varepsilon(x)$, 
for all $x,y,z\in H$.
Then we have a new Hopf algebra $H_\sigma$, the coalgebra $H$ with multiplication conjugated by $\sigma$, i.~e.
\begin{align}\label{eq:cocycle-bicomod-alg}
x \cdot_\sigma y &= 
\sigma(x_{(1)},y_{(1)}) \, x_{(2)} y_{(2)} \, \sigma^{-1}(x_{(3)},y_{(3)}), & x,y &\in H.
\end{align}

\medbreak Dually,
a twist for $H$ is an element $F = F\^{1} \otimes F\^{2} \in H\ot H$ that satisfies
\begin{align*}
(1\ot F)(\id\ot\Delta)(F)&=(F\ot 1) (\Delta\ot\id)(F), &
(\id\ot\varepsilon)(F)&=(\varepsilon\ot\id)(F)=1
\end{align*}
and has an  inverse $F^{-1} = F\^{-1} \otimes F\^{-2}$.
Then $H^F$, the algebra $H$ with the comultiplication $\Delta^F\coloneqq F\Delta F^{-1}$ is a Hopf algebra. 
These definitions are compatible with duals, 
i.e. the transpose of a twist is a cocycle etc.

\begin{example}\label{exa:twist-graded}
Let now $\varGamma$ be a group with unit $e$ and $A$ a Hopf algebra.
Assume that $A = \oplus_{\kappa \in \varGamma} A_{\kappa}$ is a $\varGamma$-graded algebra, each 
homogeneous component $A_{\kappa}$ is a subcoalgebra and $\Ss(A_{\kappa}) = A_{\kappa^{-1}}$.
Then $A_{e}$ is a Hopf subalgebra of $A$. 
Let $F = F\^{1} \otimes F\^{2} \in A_{e}\ot A_{e}$ be a twist for $A_{e}$. 
Then $F$ is a  twist for $A$, which remains a $\varGamma$-graded algebra: 
$A^F = \oplus_{\kappa \in \varGamma} A_{\kappa}^F$,
where  $A_{\kappa}^F$ is  again a sucoalgebra, namely $A_{\kappa}$ with comultiplication  conjugated by $F$.
\end{example}

\section{Hopf algebras with a central Hopf subalgebra}\label{sec:hopf-systems}
\subsection{Hopf systems}
Let  $A$ be  a central Hopf subalgebra of a Noetherian Hopf algebra $H$ and
let $G = \Alg (A, \ku)$ be the pro-affine algebraic group defined by $A$; its unit is the counit $\varepsilon$.
Given $\kappa \in G$, let
\begin{align*}
\Mg_\kappa &= \ker \kappa \in \Specmax A, & \Ig_\kappa &= H \Mg_\kappa= \Mg_\kappa H, & H_\kappa &= H/  \Ig_\kappa.
\end{align*}
Since $A$ is central, $H_{\kappa}$ is an algebra (with multiplication $m_{\kappa}$ and unit $u_{\kappa}$)
and the natural projection $p_\kappa : H \to H_\kappa$ is an algebra map.  
If $\kappa, \gamma \in G$, then 
\begin{align*}
\Delta (\Mg_{\kappa\gamma}) \subset \Mg_\kappa \otimes A + A \otimes \Mg_\gamma &\implies 
\Delta (\Ig_{\kappa\gamma}) \subset \Ig_\kappa \otimes H + H \otimes \Ig_\gamma,
\\
\Ss (\Mg_\kappa) = \Mg_{\kappa^{-1}} &\implies \Ss (\Ig_\kappa) = \Ig_{\kappa^{-1}} .
\end{align*}
Hence for any $\kappa,\gamma \in G$ there are well-defined algebra morphisms
\begin{align}\label{eq:hopf-system}
\Delta_{\kappa,\gamma}: H_{\kappa\gamma} &\to H_\kappa \otimes H_\gamma,
\\
\Ss_{\kappa}: H_{\kappa} &\to H_{\kappa^{-1}}^{\operatorname{op}},
\end{align}
that satisfy 
\begin{align}\label{eq:hopf-system-comm-diag}
&\xymatrix{H \ar@{->}[rr]^{\Delta}\ar@{->}[d]_{p_{\kappa\gamma}}  & & H \otimes H \ar@{->}[d]^{p_{\kappa} \otimes p_{\gamma}} 
\\  H_{\kappa\gamma}  \ar@{->}[rr]^{\Delta_{\kappa,\gamma}} & & H_\kappa \otimes H_\gamma,}
& &\xymatrix{H \ar@{->}[rr]^{\Ss}\ar@{->}[d]_{p_{\kappa}}  & & H^{\operatorname{op}}\ar@{->}[d]^{p_{\kappa^{-1}}}
\\  H_{\kappa}  \ar@{->}[rr]^{\Ss_{\kappa}} & & H_{\kappa^{-1}}^{\operatorname{op}}.}
\end{align}

By the coassociativity and antipode axioms, for any $\kappa,\gamma,\nu \in G$ we have
\begin{align}\label{eq:hopf-system-props}
(\Delta_{\kappa,\gamma} \otimes \id_{H_\nu})\Delta_{\kappa\gamma,\nu} &=
(\id_{H_\kappa} \otimes \Delta_{\gamma,\nu})\Delta_{\kappa,\gamma \nu}: H_{\kappa\gamma \nu} \to H_\kappa \otimes H_\gamma\otimes H_\nu,
\\
(\id_{H_\kappa} \otimes \varepsilon)\Delta_{\kappa,\varepsilon} &= \id_{H_\kappa} 
= (\varepsilon \otimes \id_{H_\kappa}) \Delta_{\varepsilon, \kappa},
\\
m_{\kappa}(\id \otimes \Ss_{\kappa^{-1}})\Delta_{\kappa,\kappa^{-1}}
&= u_{\kappa} \varepsilon = m_{\kappa}(\Ss_{\kappa^{-1}} \otimes \id)\Delta_{\kappa^{-1}, \kappa}
\end{align}
In particular 
$H_{\varepsilon}$ is a quotient Hopf algebra of $H$. In other words, $(H_{\kappa})_{\kappa \in G}$ is a Hopf system
in the sense of \cite{a-can}. By Lemma \ref{lemma:exact-sequence-hopf} \ref{item:exact-sequence-hopf-fflat} 
and Remark \ref{rem:exact-sequence-hopf-fflat} \ref{item:exact-sequence-hopf-fflat2},
we have  an exact sequence of Hopf algebras  
\[A \overset{\iota}{\hookrightarrow} H \twoheadrightarrow H_{\varepsilon}.\]

\medbreak
Given $\kappa \in G$,  $H_{\kappa}$ a $H_{\varepsilon}$-bicomodule algebra via 
\begin{align*}
\varrho_{\kappa} &\coloneqq \Delta_{\kappa, \varepsilon}: H_{\kappa} \to H_\kappa \otimes H_{\varepsilon}, &
\lambda_{\kappa} &\coloneqq \Delta_{\varepsilon, \kappa}: H_{\kappa} \to H_{\varepsilon}\otimes H_\kappa;
\\
\text{clearly } \varrho_{\kappa} p_\kappa &= (p_\kappa \otimes p_\varepsilon) \Delta_{H}, &
\lambda_{\kappa} p_\kappa &= (p_\varepsilon \otimes p_\kappa) \Delta_{H}.
\end{align*}

\begin{lemma}\label{lemma:hopf-system-cleft} \cite[Lemma 3.1]{aay}
If  $H$  is $H_{\varepsilon}$-cleft, then so is $H_{\kappa}, \forall \kappa \in G$. \qed
\end{lemma}

\begin{remark} \label{rem:H-kappa-neq0}
If  $H$  is $H_{\varepsilon}$-cleft, then $H_{\kappa} \neq 0$, i.e., $\Ig_\kappa \neq H$. 
For this, consider $\xi \in\Reg_{\varepsilon}(H, A)$ satisfying  $ \xi (ah) = a \xi(h)$, for all  $a\in A, h \in H$; see Definition--Lemma \ref{def-lemma:extension-hopf}. Then $\xi(\Ig_\kappa) = \xi(\Mg_\kappa H) \subset  \Mg_\kappa$.
Thus,  if $\Ig_\kappa = H$, then $1 \in\Ig_\kappa$ and therefore $1 = \xi(1) \in \xi(\Ig_\kappa) = \Mg_\kappa$,
a contradiction. 
\end{remark}

\subsection{Tensor categories}\label{subsec:tensor-categories}
Let $\Cr$ be the full subcategory of $\rep H$ whose objects are those where $A$ acts in a semisimple way.
Thus, if $V \in \Cr$, then
\begin{align*}
V &= \bigoplus_{\kappa \in G} V_{\kappa}, & &\text{ where } &  V_{\kappa} &= \{v\in V: z\cdot v = \kappa(z)v
\quad \forall z \in A\}.
\end{align*}
Given $\kappa \in G$, we identify $\rep H_{\kappa}$ with a  subcategory $\Cr_{\kappa}$
of $\rep H$ via restriction along $p_{\kappa}$. Then $\Cr_{\kappa}$ is indeed a
full subcategory of  $\Cr$; if $V \in \Cr$, then $V_{\kappa} \in \Cr_{\kappa}$. In other words
\begin{align}\label{eq:grading-category}
\Cr &= \bigoplus_{\kappa \in G} \Cr_{\kappa}.
\end{align}

Because of \eqref{eq:hopf-system-comm-diag}, $\Cr$ is closed under tensor products and duality.
Clearly it is is a full subcategory closed under taking subquotients, hence it is a tensor subcategory of $\rep H$
\cite[4.11.1]{EGNO}. Notice however that $\Cr$ is not closed under extensions, see the discussion 
on $\rep H$ in Subsection \ref{subsec:grading-repH}.
By \eqref{eq:hopf-system-comm-diag}, we also see that \eqref{eq:grading-category}
is a grading of tensor categories, that is
\begin{align*}
\Cr _{\kappa} \otimes \Cr_{\gamma} &\to \Cr _{\kappa\gamma}, &
(\Cr _{\kappa})^{*} = \Cr_{\kappa^{-1}}.
\end{align*}

By Remark \ref{rem:H-kappa-neq0} this grading is faithful if $H$ is $H_{\varepsilon}$-cleft.
This statement is also true when  $\dim H_\varepsilon$ is finite, as we show next.

\begin{prop}\label{prop:faithful-grading} If $\dim H_\varepsilon$ is finite,
then the  grading $\Cr  = \bigoplus_{\kappa \in G} \Cr_{\kappa}$ is faithful. 
	Hence, $H_{\kappa} \neq 0$ for all $\kappa \in \varGamma$.
\end{prop}

\pf
Let $\cA$ be the full subcategory   of  $\rep A$ whose objects $W$ are semisimple,
i.e., 
\begin{align*}
	W &= \bigoplus_{\kappa \in G} W_{\kappa}, & &\text{ where } &  W_{\kappa} &= \{w\in W: z\cdot w = \kappa(z)w,
	\quad \forall z \in A\}.
\end{align*}
We have a grading 
$\cA = \bigoplus_{\kappa \in G} \cA _{\kappa}$, where $\cA _{\kappa}$
is the full subcategory 
of $\rep A$ whose objects are those $W$ with $W = W_{\kappa}$.  
This grading is faithful since the 
one-dimensional representation supported by $\kappa$ belongs to $\cA _{\kappa}$.  
 
 The restriction functor $\rep H \to \rep A$ induces a functor $F: \Cr \to \cA$. 
 We claim that $F$ is dominant, that  is, for every object $W \in \cA$ there exists an object $V \in \Cr$ 
 such that $W$ is a subobject of $F(V)$. To see this, recall that  $H$ is faithfully flat over $A$, cf. 
 Remark \ref{rem:exact-sequence-hopf-fflat} (ii). Then the inclusion $A \hookrightarrow H$ 
 induces a monomorphism $W = A \otimes_AW \hookrightarrow H \otimes_AW$ in $\rep A$, 
 the latter being the restriction of the finite dimensional $H$-module $H \otimes_AW$, by assumption. 
 Clearly, the image of $W$ is contained in the $A$-socle $V$ of $H \otimes_AW$, which is an object of $\Cr$, 
 because $A$ is central in $H$.
The faithfulness of the grading   \eqref{eq:grading-category} follows from the fact 
that $F(\Cr_{\kappa}) \subseteq \cA _{\kappa}$, for all $\kappa \in G$.
\epf

\begin{definition}
Let $\varGamma \leq G$. Then $\Cr_{\varGamma}$ denotes the full subcategory of $\Cr$ generated by  
$\Cr_{\kappa}$, $\kappa \in \varGamma$, that is the full subcategory of $\rep H$ 
with objects where $A$ acts in a semisimple way by characters in $\varGamma$. In other words
\begin{align}\label{eq:grading-category-Gamma}
\Cr_{\varGamma} &= \bigoplus_{\kappa \in \varGamma} \Cr_{\kappa}.
\end{align}
\end{definition}

\subsection{Matrix coefficients}\label{subsec:Matrix-coefficients}
Let $R$ be an algebra.
Given $V \in \rep R$ with associated representation $\vartheta_V: R \to \End V$, 
the image of the transpose map $^{t}\vartheta_V:   (\End V)^{*} \to  R^{*}$ is denoted by $C_{V}$; 
its elements are the matrix coefficients of $V$.
Clearly $C_{V}$ is a subcoalgebra of the finite dual $R^{\circ}$ and 
\begin{align*}
R^{\circ} &= \bigcup_{V \in \rep R} C_{V}.
\end{align*}
Given a morphism of  algebras $\phi: R \to T$, the transpose  $^{t}\phi:   T^{*} \to  R^{*}$
induces a morphism of  coalgebras $\phi^{\circ}:   T^{\circ} \to  R^{\circ}$. 
Clearly, if $V \in \rep T$ and $V^{\phi} \in \rep R$ is obtained by restriction along  $\phi$, then
$\phi^{\circ}(C_V) = C_{V^{\phi}}$. 

If in addition $R$ and $T$ are  Hopf algebras  and $\phi$ is a morphism of Hopf algebras, then so is $\phi^{\circ}:   T^{\circ} \to  R^{\circ}$. 

\medbreak  
We consider the subcoalgebras of $H^{\circ}$
\begin{align}\label{eq:def-Ckappa}
C(\kappa) &\coloneqq \bigcup_{W\in \rep H_{\kappa}} C_{W^{p_{\kappa}}}, \quad \kappa \in G.
\end{align}
In other words $p_{\kappa}^{\circ}:  H_{\kappa}^{\circ} \to C(\kappa)$ is an isomorphism of coalgebras.
In particular   $C(\varepsilon) \simeq H_{\varepsilon}^{\circ}$ as Hopf algebras.

\medbreak
Let $\varGamma$ be any subgroup of $G$.  We introduce
\begin{align}\label{eq:def-HGamma}
\alg(\varGamma) \coloneqq \sum_{\kappa \in \varGamma} C(\kappa).
\end{align}

Then $\alg(\varGamma) $ is a Hopf subalgebra of $H^{\circ}$ by \eqref{eq:hopf-system-comm-diag}. 
Notice that  the subcategory $\Cr_{\kappa}$ in Subsection \ref{subsec:tensor-categories} is equivalent to the 
 category of  finite-dimensional right $C(\kappa)$-comodules.
Also $\Cr_{\varGamma}$ is tensor-equivalent to the tensor category of  finite-dimensional $\alg(\varGamma)$-comodules.  

\medbreak
Let 
$\varGamma_{\hspace{-2pt}\mathrm{fd}} \coloneqq \left\{\kappa \in \varGamma: \rep H_{\kappa}\neq 0\right\}
= \left\{\kappa \in \varGamma: C(\kappa) \neq 0\right\}$.
Clearly $\varGamma_{\hspace{-2pt}\mathrm{fd}}$ is a subgroup of  $\varGamma$ but it could be strictly smaller.

\begin{example}
Let $\mathfrak h_n$ be the $n$-th Heisenberg Lie algebra, with basis $x_i, y_i, z$, $i \in \I_{n}$ where $z$ is central,
$[x_i, x_j] = [y_i, y_j]  = 0$, $[x_i, y_j] = \delta_{ij}z$, $i, j \in \I_{n}$.  By Lemma \ref{lemma:exact-sequence-hopf}
\ref{item:exact-sequence-hopf-fflat2}, we have an exact sequence of Hopf algebras
\begin{align*}
\ku[z] \hookrightarrow U(\mathfrak h_n) \twoheadrightarrow R
\end{align*} 
where $R$ is the polynomial ring in $2n$ variables. Now $\G = \Alg (\ku[z], \ku)$ is the algebraic group $(\ku, +) $; for any 
$\kappa \in \G \setminus 0$, $H_{\kappa} =  U(\mathfrak h_n) / \langle z-\kappa\rangle$ is isomorphic to the Weyl algebra in $2n$  variables
which has no non-zero 
finite-dimensional representation, being simple and infinite-dimensional. 
Thus for any $\varGamma \leq \G$, 
$\varGamma_{\hspace{-2pt}\mathrm{fd}}$ is  trivial.
\end{example}

The following is  the main result of this Section. Recall that  $\iota: A \to H$ is  the inclusion. 
Let $\varpi: \alg(\varGamma) \to A^{\circ}$ be the restriction of 
$\iota^{\circ}: H^{\circ} \to  A^{\circ}$ and let 
$\imath: H_{\varepsilon}^{\circ} \simeq C(\varepsilon)\to \alg(\varGamma)$ be 
the  inclusion. 

\begin{theorem}\label{thm:H(Gamma)} The algebra
$\alg(\varGamma) =  \bigoplus_{\kappa \in \varGamma_{\hspace{-2pt}\mathrm{fd}}} C(\kappa)$ is faithfully  $\varGamma_{\hspace{-2pt}\mathrm{fd}}$-graded and 
the following  is an exact sequence of Hopf algebras:
\begin{align}\label{eq:exact-sequence-cokernel- group}
H_{\varepsilon}^{\circ} \xhookrightarrow[]{\imath}  \alg(\varGamma)  \xrightarrowdbl[]{\varpi} \ku \varGamma_{\hspace{-2pt}\mathrm{fd}}.
\end{align}
\end{theorem}

\pf Let $(\kappa_j)$ be a finite family of different elements in $G$. Then 
\begin{align*}
\bigcap_{j\neq i} \Mg_{\kappa_j} + \Mg_{\kappa_i} \overset{(\star)}{=}A \overset{(\#)}{\implies}
\bigcap_{j\neq i} \Ig_{\kappa_j} + \Ig_{\kappa_i} = H \overset{(\ast)}{\implies} 
\Big(\sum_{j\neq i} C(\kappa_j)\Big) \cap C(\kappa_i) = 0.
\end{align*}
Here $(\star)$ is a standard fact in commutative algebra and $(\#)$ is evident. Let us prove $(\ast)$:
If $f \in C(\kappa)$, then $f_{\vert \Ig_{\kappa}} = 0$. Thus if $f\in \Big(\sum_{j\neq i} C(\kappa_j)\Big) \cap C(\kappa_i)$,
 then $f_{\vert \bigcap_{j\neq i} \Ig_{\kappa_j} + \Ig_{\kappa_i} } = 0$, hence $f=0$.
Therefore 
 the sum $\sum_{\kappa \in \varGamma} C(\kappa)$ is direct. Since the multiplication of $H^{\circ}$ is the transpose of $\Delta$,
\eqref{eq:hopf-system-comm-diag} implies that 
$C(\kappa) \cdot C(\gamma) \subseteq  C(\kappa\gamma)$, for all  $\kappa, \gamma \in \varGamma$.

\medbreak
We next claim that $\Imm \varpi \simeq \ku \varGamma_{\hspace{-2pt}\mathrm{fd}}$, so that the map $\varpi$ in \eqref{eq:exact-sequence-cokernel- group} makes sense. 
Indeed, if $V\in \rep H_{\kappa}$, 
then  $A$ acts on $V$ via $\kappa$ and thus $C_{\Res_A V} \subseteq\ku \kappa$ 
and the equality holds iff $V\neq 0$. This implies that
$\Imm \varpi \subseteq \ku \varGamma_{\hspace{-2pt}\mathrm{fd}}$ and the equality follows by definition of $\varGamma_{\hspace{-2pt}\mathrm{fd}}$.

\medbreak  
Our next goal is to show that 
\begin{align}\label{eq:exact-sequence-cokernel- group-coinv}
H_{\varepsilon}^{\circ} = \alg(\varGamma)^{\operatorname{co} \varpi}.
\end{align}

We start with the following observation.
By standard arguments on matrix coefficients, if
$W \in \rep H_{\kappa}$ and $V = W^{p_{\kappa}}$, then
\begin{align}\label{eq:counit}
\varpi (f)  &= \boldsymbol{\varepsilon}(f) \kappa, & f &\in C_{V}
\end{align}
where the counit of $H^{\circ}$ is denoted by $\boldsymbol{\varepsilon}$.
Now \eqref{eq:exact-sequence-cokernel- group-coinv} follows because each $C(\kappa)$ is a subcoalgebra: 
given $f\in \alg(\varGamma)$, write $f =  \sum_{\kappa \in \varGamma} f_{\kappa}$
with $f_{\kappa} \in C(\kappa)$. Then by \eqref{eq:counit},
\begin{align*}
(\id \otimes \varpi) \Delta (f) =   \sum_{\kappa \in \varGamma} f_{\kappa} \otimes \kappa,
\end{align*}
hence $f \in \alg(\varGamma)^{\operatorname{co} \varpi}$ if and only if $f \in C(\varepsilon)$.

\medbreak
For the exactness of  \eqref{eq:exact-sequence-cokernel- group}, it remains to see that
$\ker\varpi = \alg(\varGamma)  (H_{\varepsilon}^{\circ})^+$. 
Now $\alg(\varGamma)$ is faithfully coflat over $\ku \varGamma_{\hspace{-2pt}\mathrm{fd}}$ by Lemma \ref{lema:exact-sequence-hopf-fcoflat}.
By Lemma \ref{lemma:exact-sequence-hopf2} \ref{item:exact-sequence-hopf-fcoflat} we are reduced to show that
$\varpi$ is conormal, but this follows from the centrality of $A$: take $f\in \alg(\varGamma)$, then 
$(\id \otimes \varpi)\varrho_{\ell}(f) = f\_{1} \Ss (f\_{3}) \otimes \varpi (f\_{2}) = 0 $
if and only if   $\langle f\_{1} \Ss (f\_{3}), h\rangle \langle \varpi (f\_{2}), a\rangle =0$ for all $h\in H$, $a \in A$.
Now
\begin{align*}
\langle f\_{1} \Ss (f\_{3}), h\rangle \langle \varpi (f\_{2}), a\rangle  &=
\langle f\_{1}, h\_{1}\rangle  \langle f\_{3}, \Ss (h\_{2})\rangle \langle \varpi (f\_{2}), a\rangle 
\\ &= \langle f, h\_{1} a \Ss (h\_{2})\rangle = \varepsilon(h)  \langle f,a\rangle.
\end{align*}
Thus, if $f \in \ker \varpi$,  then  $f\_{1} \Ss (f\_{3}) \otimes \varpi (f\_{2}) = 0 $, hence $\ker \varpi$ is a
subcomodule with respect to $\varrho_{\ell}$, i.e. $\varpi$ is conormal.
\epf

\begin{remark}\label{rem:varGamma-subgroup}
If $\varGamma' \leq \varGamma$, then  $\alg(\varGamma')$ is a Hopf subalgebra of  $\alg(\varGamma)$,
$\varGamma'_{\hspace{-2pt}\mathrm{fd}} \leq \varGamma_{\hspace{-2pt}\mathrm{fd}}$ 
and there is a morphism of exact sequences
\begin{align*}
\xymatrix{ H_{\varepsilon}^{\circ} \ar@{=}[d] \ar@{^{(}->}[r]  & \alg(\varGamma') \ar@{->>}[r]  \ar@{^{(}->}[d]
& \ku \varGamma'_{\hspace{-1pt}\mathrm{fd}} \ar@{^{(}->}[d]
\\
H_{\varepsilon}^{\circ}  \ar@{^{(}->}[r]  & \alg(\varGamma) \ar@{->>}[r]  & \ku \varGamma_{\hspace{-1pt}\mathrm{fd}}.}
\end{align*}
\end{remark}

\begin{remark}\label{rem:cocentral} The projection $\alg(\varGamma) \to \ku \varGamma_{\hspace{-1pt}\mathrm{fd}}$ is cocentral.
\end{remark}

\pf This follows from the centrality of $A$ in $H$.
\epf

Finally we establish some properties of $\alg(\varGamma)$, see \S \ref{subsec:cofrob}. First we recall:

\begin{theorem}\label{thm:exseq2} \cite[2.10, 2.13]{AC} 
Let $A\xhookrightarrow[]{\iota} C \xrightarrowdbl[]{\pi} B$ be an exact sequence of Hopf algebras with $C$ faithfully coflat as a $B$-comodule. Then
$C$ is co-Frobenius (respectively cosemisimple) if and only if $A$ and $B$ are co-Frobenius (respectively cosemisimple).
\qed
\end{theorem}

Theorem \ref{thm:exseq2} applied to the exact sequence \eqref{eq:exact-sequence-cokernel- group},
see Lemma \ref{lema:exact-sequence-hopf-fcoflat}, gives:

\begin{coro}\label{coro:exseq2}
$\alg(\varGamma)$ is co-Frobenius (respectively cosemisimple) if and only if $H_{\varepsilon}^{\circ}$ is co-Frobenius
(respectively cosemisimple). \qed
\end{coro}

\subsection{A grading of $\rep H$}\label{subsec:grading-repH}
Given $V \in \rep H$, we set
\begin{align*}
V_{(\kappa)} &= \{v\in V: (z - \kappa(z))^n\cdot v = 0\quad \forall z \in A, \quad \text{for some } n \in \N\}.
\end{align*}
Hence $V = \bigoplus_{\kappa \in G} V_{(\kappa)}$ by a classical argument.
Given $\kappa \in G$, we introduce the full subcategory $\Cr_{(\kappa)}$
of $\rep H$ whose objects are those $V \in \rep H$ with $V = V_{(\kappa)}$.   Therefore
\begin{align}\label{eq:grading-repH}
\rep H &= \bigoplus_{\kappa \in G} \Cr_{(\kappa)}.
\end{align}

Since $\Cr_{\kappa} \subset \Cr_{(\kappa)}$, the grading \eqref{eq:grading-category} of $\Cr$
is inherited from \eqref{eq:grading-repH}.

\begin{prop}\label{prop:grading-repH}  The decomposition \eqref{eq:grading-repH} is a grading of tensor categories.
If $\dim H_{\varepsilon} < \infty$, then this grading  is faithful.
\end{prop}

\pf  Let $V, W\in \rep H$. We first observe that $V \in \Cr _{(\kappa)}$ iff there exists a filtration
\begin{align*}
0 = V_{-1}\subset V_{0} \subset V_{1} \subset  \dots \subset V_{n} = V
\end{align*}
such that $V_{i}/ V_{i-1} \in \Cr_{\kappa}$ for all $i\in \I_{0,n}$. Thus, if $V \in \Cr _{(\kappa)}$ 
and $W \in \Cr _{(\gamma)}$  with an analogous filtration 
$0 = W_{-1}\subset W_{0} \subset W_{1} \subset  \dots \subset W_{m} = W$, then 
the filtration  of $V \otimes W$  given by $(V\otimes W)_{-1} = 0$ and 
\begin{align*}
(V \otimes W)_{j} &= \sum_{r,s\geq 0, r+s = j} V_{r}\otimes W_{s}, \quad j \geq 0,
\end{align*}
is exhaustive and satisfies $(V \otimes W)_{j}/ (V \otimes W)_{j-1} \in \Cr_{\kappa\gamma}$.
Indeed, given $r$ and $s$, take $v \in V_{r}$ and $w\in W_{s}$. 
By hypothesis, for any $z \in A$ we have
\begin{align*}
z\cdot v &\in \kappa(z) v + V_{r-1}, & z\cdot w &\in\gamma(z) w + W_{s-1}, 
\end{align*}
hence
\begin{align*}
z\cdot (v\otimes w) &= z\_{1}\cdot v  \otimes z\_{2}\cdot w  \in 
 (\kappa(z\_{1}) v + V_{r-1}) \otimes (\gamma(z\_{2}) w + W_{s-1}) 
 \\ &\subseteq  \kappa\gamma(z) (v\otimes w) + V_{r}\otimes W_{s-1} + V_{r-1}\otimes W_{s}
 + V_{r-1}\otimes W_{s-1}
 \\ &\subseteq  \kappa\gamma(z) (v\otimes w) + (V\otimes W)_{r+s-1}.
\end{align*}
Thus
$\Cr _{(\kappa)} \otimes \Cr_{(\gamma)} \hookrightarrow \Cr _{(\kappa\gamma)}$.
Now, if $v\in V_{i}$ and $f\in V^*$, then for any $z \in A$ we have
$\langle z\cdot f, v\rangle = \langle f, \Ss(z)\cdot v\rangle = \langle f, \kappa^{-1}(z)\cdot v + v'\rangle $
for some $v'\in V_{i-1}$. Hence 
\begin{align*}
\langle (z - \kappa^{-1}(z))\cdot f, v\rangle = \langle f, v'\rangle.
\end{align*}
Iterating we see that $(z - \kappa^{-1}(z))^{n + 1}\cdot f = 0$, for $n$ as above. Therefore
$(\Cr _{(\kappa)})^{*} = \Cr_{(\kappa^{-1})}$.

\medbreak
Observe that the previous arguments imply that the category $\rep A$ also bears a grading 
$\rep A = \bigoplus_{\kappa \in G} \cA _{(\kappa)}$, where $\cA _{(\kappa)}$
is the full subcategory 
of $\rep A$ whose objects are those $V$ with $V = V_{(\kappa)}$.  
Now this grading is faithful since the 
one-dimensional representation supported by $\kappa$ belongs to $\cA _{(\kappa)}$.  

\medbreak 
Assume now that $\dim H_{\varepsilon} < \infty$. As in Proposition \ref{prop:faithful-grading}, the restriction functor $F: \rep H \to \rep A$ is dominant, implying the faithfulness of the grading   \eqref{eq:grading-repH}.
\epf

\begin{definition}
Let $\varGamma \leq G$. Then $\Cr_{(\varGamma)}$ denotes the full subcategory of $\rep H$ generated by  
$\Cr_{(\kappa)}$, $\kappa \in \varGamma$, in other words
$\Cr_{(\varGamma)} = \bigoplus_{\kappa \in \varGamma} \Cr_{(\kappa)}$.
We introduce
\begin{align*}
\alg((\varGamma)) \coloneqq \sum_{W\in \Cr_{(\varGamma)}} C_{W};
\end{align*}
this is a Hopf subalgebra of $H^{\circ}$ and
 $\Cr_{(\varGamma)}$ is tensor-equivalent to the tensor category of  finite-dimensional $\alg((\varGamma))$-comodules.  
\end{definition}

The Hopf algebra $\alg((\varGamma))$ will be studied elsewhere.

\section{Finite-by-cocommutative Hopf algebras}\label{sec:finite-by-coco}

\subsection{Discussion of the assumptions}\label{subsec:assumptions}
As in the previous Section we fix a Noetherian Hopf algebra $H$ with a central Hopf subalgebra $A$; 
thus  we have  an exact sequence of Hopf algebras  
$(\E)$: $A \overset{\iota}{\hookrightarrow} H \overset{p_{\varepsilon}}{\twoheadrightarrow} H_{\varepsilon}$.
Let $G = \Alg(A, \ku)$.

\subsubsection{Assumptions on $H_{\varepsilon}$}

Consider the following assumptions:

\begin{assumption}\label{assumption:fingen}
$H$ is a finitely generated $A$-module, say $(h_i)_{i\in \I_N}$ generate $H$ over $A$. 
\end{assumption}

\begin{assumption}\label{assumption:findim}
$H_{\varepsilon}$ is finite-dimensional. 
\end{assumption}

\begin{remark} Assumption \ref{assumption:fingen} implies \ref{assumption:findim}. In general, we have
for any $\kappa \in G$
\begin{align}\label{eq:dimension-bound}
\dim H_{\kappa} & \leq N.
\end{align}
Conversely, Assumption \ref{assumption:findim} implies Assumption \ref{assumption:fingen} if moreover the extension 
$(\E)$ is \emph{cleft}.
\end{remark}

\begin{question}\label{question:finiteness-assumptions}
Does  Assumption \ref{assumption:findim} imply Assumption \ref{assumption:fingen} always?
\end{question}

\begin{example}
Assumption   \ref{assumption:findim} does not imply cleftness of $(\E)$; see the example by Oberst and Schneider in 
\cite[Section 3]{montgomery}. 
\end{example}

If \eqref{eq:dimension-bound} holds, then the subcoalgebra $C(\kappa)$ of the finite dual
$H^{\circ}$ defined in \eqref{eq:def-Ckappa} is identified with the finite-dimensional coalgebra $H_{\kappa}^{*}$, 
via the injective map of coalgebras $p_{\kappa}^t : H_{\kappa}^{*} \to H^{\circ}$.
In particular,  $C(\varepsilon) \simeq H_{\varepsilon}^{*} = H_{\varepsilon}^{\circ}$.

\begin{lemma}If Assumption \ref{assumption:findim} holds, then  $\dim H_{\kappa} = \dim H_{\varepsilon}$, for all $\kappa \in G$.
\end{lemma}

Compare with \cite[Theorem 3.5.2]{EGNO}.

\pf
First recall that for any Hopf algebra $K$ and  right $K$-comodule $V$, one has an isomorphism of right $K$-comodules
\begin{align*}
V\otimes K &\to V_{\text{trivial}} \otimes K, & v \otimes k &\mapsto v\_0 \otimes v\_{1} k, &
v\in V, \ k &\in K. 
\end{align*}
Hence, if $\dim V < \infty$, then $V\otimes K \simeq K^{\dim V}$. Similarly, $K\otimes V  \simeq K^{\dim V}$.
Consider the coalgebra decomposition 
\[\alg(\varGamma) = \bigoplus_{\kappa \in \varGamma} C(\kappa),\]
where $\dim C(\kappa) < \infty$ by \eqref{eq:dimension-bound}.
Recall that $\Cr_{\kappa}$ is  the category of finite-dimensional $C(\kappa)$-comodules. 
By the previous discussion, for any $\gamma \in \varGamma$,
\begin{align*}
\bigoplus_{\kappa \in \varGamma} \left(C(\kappa) \otimes C(\gamma) \right)
\simeq \alg(\varGamma)\otimes C(\gamma)  \simeq \alg(\varGamma)^{\dim C(\gamma)}   
\simeq \bigoplus_{\kappa \in \varGamma} C(\kappa) ^{\dim C(\gamma)}  
\end{align*}
as $\alg(\varGamma)$-comodules. Notice that $C(\kappa) \otimes C(\gamma) \in \Cr_{\kappa\gamma}$;
while $C(\kappa) ^{\dim C(\gamma)} \in \Cr_{\kappa}$. Therefore 
\begin{align}\label{eq:C-kappa-gamma}
C(\kappa) \otimes C(\gamma) &\simeq C(\kappa\gamma) ^{\dim C(\gamma)},&
\text{for all } \kappa, \gamma \in \varGamma.
\end{align}
By a similar argument, using the isomorphism $C(\kappa) \otimes \alg(\varGamma)\simeq \alg(\varGamma)^{\dim C(\kappa)}$, we also obtain that 
\begin{align}\label{eq:C-gamma-kappa}
	C(\kappa) \otimes C(\gamma) &\simeq C(\kappa\gamma) ^{\dim C(\kappa)},&
	\text{for all } \kappa, \gamma \in \varGamma.
\end{align}
Now $\dim C(\kappa) \neq 0$, by Proposition \ref{prop:faithful-grading}.
Comparing dimensions on  \eqref{eq:C-kappa-gamma} and \eqref{eq:C-gamma-kappa}, we get $\dim C(\kappa) = \dim C(\gamma)$, for all $\kappa, \gamma  \in \varGamma$. This proves the lemma. 
\epf

\emph{In the rest of the paper, we shall assume that Assumption \ref{assumption:findim} holds and that the extension 
$(\E)$ is cleft, hence Assumption \ref{assumption:fingen} also holds. }

\subsubsection{Assumptions on $H$} Let $H$ be a Hopf algebra  with a central Hopf subalgebra $A$
such that $H$ is a finitely generated $A$-module.  Then 
\begin{align*}
H \text{ is Noetherian} \iff A \text{ is Noetherian} \iff A \text{ is affine.} 
\end{align*}
 In particular, $G = \Alg(A, \ku)$ is an algebraic group.
These equivalences follow from Lemma \ref{lemma:A-affine} below whose proof requires two results.

\begin{theorem}\label{thm:formanek-jategaonkar}\cite{formanek-jategaonkar}
If $R$ is a right Notherian ring which is finitely generated as a right module
over a commutative subring $S$, then $S$ is Noetherian. \qed
\end{theorem}

\begin{theorem}\label{thm:molnar}\cite{molnar}
Let  $A$ be a commutative Hopf algebra. Then $A$ is Noetherian if and only if it is affine. \qed
\end{theorem}

\begin{lemma}\label{lemma:A-affine}  
Let
$H$ be a Hopf algebra  with a commutative Hopf subalgebra $A$
such that $H$ is a finitely generated $A$-module.
The following are equivalent:
\begin{enumerate}[leftmargin=*,label=\rm{(\alph*)}] 
\item\label{item:H-noeth} $H$ is Noetherian.

\item\label{item:A-noeth} $A$ is Noetherian.

\item\label{item:A-affine} $A$ is affine.
\end{enumerate}
\end{lemma}

\pf \ref{item:H-noeth} $\Rightarrow$ \ref{item:A-noeth} is Theorem \ref{thm:formanek-jategaonkar}.
\ref{item:A-noeth} $\Rightarrow$ \ref{item:H-noeth}: By assumption $H$ is a Noetherian $A$-module 
and so it is a Noetherian algebra. 
\ref{item:A-noeth} $\Leftrightarrow$ \ref{item:A-affine} is Theorem \ref{thm:molnar}.
\epf

\subsubsection{Structure of $\alg(\varGamma)$} 

Let $\varGamma \leq G$; then $\varGamma_{\hspace{-2pt}\mathrm{fd}} = \varGamma$ because of \ref{eq:dimension-bound}  and Remark \ref{rem:H-kappa-neq0}.
Our main object of interest  is the Hopf algebra
$\alg(\varGamma) =  \bigoplus_{\kappa \in \varGamma} C(\kappa)$. By Theorem \ref{thm:H(Gamma)}  
it fits into the exact sequence 
\begin{align}
\tag{$\F^{\varGamma}$} H_{\varepsilon}^{*} \xhookrightarrow[]{\imath}  \alg(\varGamma)  \xrightarrowdbl[]{\varpi} \ku \varGamma.
\end{align}

\begin{lemma}\label{lemma:Gamma-cleft} The extension $(\F^{\varGamma})$ is cleft; hence 
$\alg(\varGamma)$ is strongly $\varGamma$-graded.
\end{lemma}
\pf
Let $\chi: H_{\varepsilon} \to H$ be a cleaving map and let $\boldsymbol{\xi}: \alg(\varGamma) \to  H_{\varepsilon}^*$ 
be the restriction of the transpose $\chi^{t}: H^* \to H_{\varepsilon}^*$. 
Then $\boldsymbol{\xi}$ satisfies the conditions of Definition--Lemma \ref{def-lemma:extension-hopf} \ref{item:xi}, hence
$(\F^{\varGamma})$ is cleft. The second claim follows from Theorem \ref{thm:doi-takeuchi} and Example \ref{exa:strongly-graded}.
\epf

\begin{remark}\label{rem:H(Gamma)affine}
If the family $(\kappa_i)_{i \in I}$ generates $\varGamma$, then 
$\bigoplus_{i \in I} C(\kappa_i)$ generates the algebra $\alg(\varGamma)$. Hence, 
if $\varGamma$ is finitely generated, then   $\alg(\varGamma)$ is affine.
\end{remark}

\subsubsection{Pairing} 
The natural pairing $\langle\,,\,\rangle: H^* \times H \to \ku$ descends to a 
pairing $\langle\,,\,\rangle: \alg(\varGamma) \times H \to \ku$. Clearly  it is non-degenerate on one side.
Let
$I = \alg(\varGamma)^{\perp} \subseteq H$;
this is a Hopf ideal since $\langle\,,\,\rangle$ is a Hopf pairing. We have
\begin{align*}
I =  \cap_{\kappa \in \varGamma} C(\kappa)^{\perp} =  \cap_{\kappa \in \varGamma} \Ig_\kappa,
\end{align*}
by \eqref{eq:def-HGamma} and since all $H_\kappa$'s are finite-dimensional. Let $\kappa \in \varGamma$,  $W\in \rep H_{\kappa}$ and let $\ct \in C_{W^{p_{\kappa}}}$ be a matrix coefficient.
Then
$\langle\ct, a\rangle = \kappa(a) \ct$,   for any $a \in A$.
Hence 
\begin{align*}
I  \cap A=  \cap_{\kappa \in \varGamma} C(\kappa)^{\perp}\cap A 
=   \left\{ a \in A:  \kappa(a) = 0 \text{ for all } \kappa \in \varGamma \right\}.
\end{align*}
That is, $I \cap A$ is the ideal of functions on $G$ that vanish in the Zariski closure $\overline{\varGamma}$ of $\varGamma$. 
Hence, if $\overline{\varGamma} \neq G$, then $ I \cap A \neq 0$ and so $I \neq 0$. We have proved:

\begin{lemma}
If the pairing $\langle\,,\,\rangle: \alg(\varGamma) \times H \to \ku$ is non-degenerate, then 
$\varGamma$ is Zariski-dense in $G$.  \qed
\end{lemma}

\begin{question} If $\varGamma$ is Zariski-dense in $G$, is 
the pairing $\langle\,,\,\rangle: \alg(\varGamma) \times H \to \ku$  non-degenerate?
That is, does $ I \cap A = 0$ imply that $I = 0$?
\end{question}

\subsubsection{Comparison} Consider another cleft exact sequence of Hopf algebras  
$A' \overset{\iota}{\hookrightarrow} H' \overset{p_{\varepsilon}}{\twoheadrightarrow} H'_{\varepsilon}$ 
with $H^{\prime }$ Noetherian and $\dim H'_{\varepsilon} < \infty$; set $G' = \Alg(A', \ku)$. Given subgroups $\varGamma \leq G$ and 
$\varGamma' \leq G'$ it is natural to ask whether the Hopf algebras $\alg(\varGamma)$ and  $\alg(\varGamma')$ are isomorphic
(here is an abuse of notation, as the initial data are different).
We offer a partial answer, applicable in many examples.
Recall that $\varpi: \alg(\varGamma)\ku \varGamma$ is cocentral, cf. Remark \ref{rem:cocentral}.

\begin{prop}\label{prop: hopf-center} Let $f: \alg(\varGamma)\to \alg(\varGamma')$ be an isomorphism of Hopf algebras. 
If  $\hcocenter(\alg(\varGamma))= \ku \varGamma$ and $\hcocenter(\alg(\varGamma')) = \ku \varGamma'$, then $f$ induces isomorphisms  $\varphi: \varGamma \to\varGamma'$ and  $f_0: H^*_{\varepsilon}\simeq {H_{\varepsilon}^{\prime }}^*$ such that the following diagram commutes:
\begin{align*}
	\xymatrix{ H^*_{\varepsilon} \ar@{->}[d]^{f_0}  \ar@{^{(}->}[r]  & \alg(\varGamma) \ar@{->>}[r]^{\varpi}  \ar@{->}[d]^{f}
		& \ku \varGamma \ar@{->}[d]^{\varphi}
		\\
		{H_{\varepsilon}^{\prime}}^*  \ar@{^{(}->}[r]  & \alg(\varGamma') \ar@{->>}[r]^{\varpi'}  & \ku \varGamma'.}
\end{align*}
\end{prop}

\pf
Since the morphism $\varpi' f: \alg(\varGamma) \to \ku \varGamma'$ is cocentral, there exists a unique morphism of Hopf algebras  $\varphi: \ku \varGamma\to \ku\varGamma'$ such that $\varphi \varpi = \varpi' f$.
Then $\varphi$ must be an isomorphism which, combinded with the fact that $H^*_{\varepsilon} = \alg(\Gamma)^{co \varpi}$ and ${H_{\varepsilon}^{\prime}}^* = \alg(\Gamma')^{co \varpi'}$, implies that  $f$ induces by restriction an isomorphism $f_0: H^*_{\varepsilon}\simeq {H_{\varepsilon}^{\prime }}^*$ making the diagram commute. 
\epf

\begin{remark} By  Lemma \ref{lema:trivial-center} \ref{item:trivial-cocenter},  the assumptions of Proposition \ref{prop: hopf-center} are satisfied if $H_\varepsilon$ and $H_{\varepsilon}^{\prime}$ are simple and noncommutative or, more generally, if  $\hcenter(H_{\varepsilon})\simeq \ku \simeq \hcenter(H_{\varepsilon}^{\prime })$.
\end{remark}

\subsection{Coradical filtration and co-Frobenius property}\label{subsec:cofrob}
Recall that a Hopf algebra $K$ is co-Frobe\-nius if it admts a non-zero (right) integral, i.e., a linear functional
$\int: K \to \ku$, $\int \neq 0$, invariant under the dual of the left regular representation. See e.g.
\cite{AC,ace} for details and a list of equivalent characterizations. 

Let $\varGamma \leq G$. 
Corollary \ref{coro:exseq2} implies that $\alg(\varGamma)$ is co-Frobenius. 
The next Theorem gives a refinement of this fact.
\begin{theorem}\label{thm:cofrob} 
The coradical filtration of $\alg(\varGamma)$ is 
\begin{align}\label{eq:corad-filt-AGamma}
\corad_n \alg(\varGamma) &=  \bigoplus_{\kappa \in \varGamma} \corad_n C(\kappa),
\end{align}
hence $\alg(\varGamma)$ is co-Frobenius.
\end{theorem}

\pf The proof of \eqref{eq:corad-filt-AGamma} is by induction on $n \geq 0$. For $n=0$ this 
is a consequence of \cite[3.4.3]{rad-libro}, while the inductive step follows from \cite[2.4.3]{rad-libro}.
Then \eqref{eq:dimension-bound} implies that $\corad_N C(\kappa)= C(\kappa)$ and therefore $\corad_N \alg(\varGamma) = \alg(\varGamma)$. Finally, by finiteness of the coradical filtration, $\alg(\varGamma)$ is co-Frobenius cf. \cite[Theorem 2.1]{andrusdasca}.
\epf

By the Larson-Radford theorem \cite{larson-radford}, Corollary \ref{coro:exseq2}  also implies the following result.

\begin{theorem}\label{thm:coss} 
$\alg(\varGamma)$ is cosemisimple if and only if $H_{\varepsilon}$ is 
 semisimple.   \qed
\end{theorem}

\subsection{Gelfand-Kirillov dimension}\label{subsec:group-finitegk}
We refer to \cite{krause-lenagan} for the definition and basic properties of this notion.

The main result of this subsection characterizes the algebras $\alg(\varGamma)$ with finite 
Gelfand-Kirillov dimension. We start by a general remark.
Let $\varGamma$ be a group 
and let $R = \bigoplus_{\kappa \in \varGamma} R_{\kappa}$ be a strongly $\varGamma$-graded algebra.

\begin{lemma} \label{lema:GK-dim-Gamma-graded}
If $\dim R_{\kappa} < \infty$ for all $\kappa \in \varGamma$,
then $\GK R \geq \GK \ku \varGamma$. Furthermore the equality holds if  there exists $N \in \N$ such that 
\begin{align}\label{eq:GK-dim-Gamma-graded}
\dim R_{\kappa} &\leq N, & \text{ for all } \kappa &\in \varGamma.
\end{align}
\end{lemma}

In particular, if $\GK R < \infty$, then $\varGamma$ has polynomial growth.

\pf First notice that $R_{\kappa} \neq 0$ for all $\kappa \in \varGamma$ since 
$0 \neq R_e = R_{\kappa}R_{\kappa^{-1}}$.
Given a finite subset $X \subset \varGamma$, set $R_{X} = \bigoplus_{\kappa \in X} R_{\kappa}$. 
Then $\vert  X\vert \leq \dim R_{X}  < \infty$ by the preceding. We have 
$R_{X}R_{Y} = R_{XY}$ for any $Y \subset \varGamma$, since $R$ is strongly graded.
Hence $(R_{X})^n = R_{X^n}$ and 
\begin{align*}
\vert  X^n \vert \overset{\spadesuit}{\leq} \dim (R_{X})^n & \overset{\clubsuit}{\leq} N \vert  X^n\vert ,& n &\in \N,
\end{align*}
where in the second inequality we assume \eqref{eq:GK-dim-Gamma-graded}. 
Now   $\spadesuit$ implies that
\begin{align*}
\limsup _{{n\to \infty }}\log _{n}	\vert  X^n \vert & \leq \GK R, & \text{ hence } \GK \ku \varGamma & \leq \GK R.
\end{align*}
Assume that   \eqref{eq:GK-dim-Gamma-graded} holds. Let $V$ be a finite-dimensional subspace of $R$;
clearly there exists a finite $X \subset \varGamma$ such that $V \subset R_{X}$. 
Then $\clubsuit$ implies that
\begin{align*}
\limsup _{{n\to \infty }}\log _{n} \dim	V^{n} \leq \limsup _{{n\to \infty }}\log _{n}\dim (R_{X})^n 
& \leq \GK \ku \varGamma, 
\end{align*}
and the equality $\GK R = \GK \ku \varGamma$ follows.
\epf

Recall that a finitely generated
group is \emph{nilpotent-by-finite}  if it has a normal nilpotent subgroup of finite index.
Here is a celebrated result by Gromov:

\begin{theorem}\label{thm:gromov} \cite{gromov} If  $\varGamma$ is a finitely generated group, then
$\GK \ku \varGamma < \infty$  if and only if  $\varGamma$  is nilpotent-by-finite. \qed
\end{theorem}

Assume next that  $\varGamma$ is a not necessarily finitely generated group. Then $\ku \varGamma$ has finite $\GK$ 
if and only if there exists $N \in \N$ such that 
\begin{align*}
	\GK \ku \varUpsilon < N    \text{ for any finitely generated } \varUpsilon \leq \varGamma.
\end{align*}
In particular any finitely generated subgroup of $\varGamma$ should be  nilpotent-by-finite.
We then conclude:

\begin{theorem}\label{thm:gkdim} Let  $\varGamma \leq G$. Then
\begin{align}
\label{eq:gkdim-thm}\GK \alg(\varGamma)  = \GK \ku \varGamma.
\end{align}
Thus $\GK \alg(\varGamma) < \infty$  
iff   there exists $N \in \N$ such that  any finitely generated  $\varUpsilon \leq \varGamma$  is nilpotent-by-finite
and $\GK \ku \varUpsilon < N$.
\end{theorem}

\pf Assume first that $\varGamma$ is  finitely generated.
Then \eqref{eq:gkdim-thm} follows from Lemma \ref{lema:GK-dim-Gamma-graded} 
since \eqref{eq:dimension-bound} gives \eqref{eq:GK-dim-Gamma-graded}. Therefore
$\GK \alg(\varGamma) < \infty$ iff $\varGamma$ has polynomial growth iff 
 $\varGamma$  is nilpotent-by-finite by Theorem \ref{thm:gromov}. 

In general any finitely generated subalgebra of  $\alg(\varGamma)$ is contained in $\alg(\varUpsilon)$ for some 
finitely generated $\varUpsilon \leq \varGamma$, cf. Remark \ref{rem:varGamma-subgroup}. Thus
\begin{align*}
\GK \alg(\varGamma)  = \sup_{\substack{\varUpsilon \leq \varGamma \\ \text{fin. gen.}}}
\GK \alg(\varUpsilon) = \sup_{\substack{\varUpsilon \leq \varGamma \\ \text{fin. gen.}}}
\GK \ku \varUpsilon = \GK \ku \varGamma.
\end{align*}
So, \eqref{eq:gkdim-thm} holds in general and the theorem follows from Theorem \ref{thm:gromov} by the considerations above.
\epf

\medbreak
We list a few examples for illustration:

\begin{enumerate}[leftmargin=*,label=\rm{(\alph*)}] 
\item $\GK \ku \mathbb Q = 1$ because any finitely generated subgroup  of $\mathbb Q$ is cyclic (and torsion-free). 
Since $\mathbb Q$ is a subgroup of the additive group $\ku$, it embeds in any algebraic group which is not a torus.

\medbreak
\item Let $\mathbb G_{\infty} \leq \kut$ be the group of all roots of 1. Then
$\GK \ku \mathbb G_{\infty} = 0$ because any finitely generated subgroup  is cyclic (and torsion). 
Now $\mathbb G_{\infty}$ embeds in any algebraic group that contains a torus.
\end{enumerate}

\subsection{Noetherianity}\label{subsec:group-Noetherian}
A general reference for this Subsection is \cite{mcconell-robson}.
Recall that a solvable group is polycyclic if every subgroup is finitely generated; equivalently, if  it admits a subnormal series with cyclic factors. 
Also,  a group is \emph{polycyclic-by-finite} if it has  a normal polycyclic subgroup of finite index. 
The Hirsch number of a polycyclic group is the number of infinite factors in any subnormal series;
the Hirsch number of a polycyclic-by-finite group is that of a polycyclic normal subgroup with finite index. 

It is a classical result that the group algebra of a polycyclic-by-finite group is Noetherian \cite{hall};
a well-known open question is whether the converse holds.  

Recall that a group is Noetherian if it satisfies the maximal condition on subgroups; if the group algebra of a given group is Noetherian, then so is the group but the converse is not true, see \cite{ivanov}. However for linear groups there is a remarkable result of Tits:

\begin{theorem}\label{thm:linear-noeth} \cite{tits}
A linear Noetherian group  is polycyclic-by-finite. \qed
\end{theorem}

In consequence, if $\varGamma$ is a linear group, then the following are equivalent:
\begin{align*}
&\ku\varGamma \text{ is Noetherian}& &\Leftrightarrow& 
&\varGamma \text{ is Noetherian} & &\Leftrightarrow& 
&\varGamma \text{ is polycyclic-by-finite.}
\end{align*}

\begin{theorem}\label{thm:nast-fvo} \cite[Thm. 5.5]{nastasescu-fvo}
If $\varGamma$ is a polycyclic-by-finite group and $R$ is a strongly
$\varGamma$-graded ring, then $R$ is right  Noetherian when $R_e$ is so. \qed
\end{theorem}

\begin{theorem}\label{thm:Noetherian}  Let $\varGamma \leq G$.
The following are equivalent:
\begin{enumerate}[leftmargin=*,label=\rm{(\alph*)}] 
\item\label{item:Noetherian-1} The algebra $\alg(\varGamma)$ is Noetherian.

\medbreak
\item\label{item:Noetherian-4} The group $\varGamma$ is polycyclic-by-finite (in particular it is solvable).
\end{enumerate}

\end{theorem}

\pf \ref{item:Noetherian-1} $\Rightarrow$ \ref{item:Noetherian-4}:
Recall that $\varGamma$, being a subgroup of an affine algebraic group, is linear.
Since $\ku \varGamma$ is Noetherian, being a  quotient of $\alg(\varGamma)$, Theorem \ref{thm:linear-noeth} applies.
\ref{item:Noetherian-4} $\Rightarrow$ \ref{item:Noetherian-1}: By Lemma \ref{lemma:Gamma-cleft}, Theorem \ref{thm:nast-fvo} and \eqref{eq:dimension-bound}.
\epf

\subsection{Regularity}\label{subsec:group-regular}
A reference for this Subsection is \cite[Chapter 7]{mcconell-robson}. As in \emph{loc. cit.} we use the following abbreviations:
 $\pd$ stands for projective dimension;  $\rgldim$, $\lgldim$, $\gldim$ stand for
right global dimension, respectively left global dimension, and  global dimension.
Given an algebra $R$, if $\rgldim R= \lgldim R$ then we set $\gldim R = \lgldim R$;
we say that $R$ is regular if $\gldim R < \infty$. 

To start with, recall from \cite{lorenz-lorenz} the estimate 
\begin{align}\label{eq:regular-gral}
\lgldim C \leq\rgldim B + \lgldim A
\end{align}
for a crossed product $C = A \#_{\sigma} B$ of an associative algebra $A$ with a Hopf algebra $B$ (with antipode not necessarily bijective). 
See \cite{dasca-nasta-torrecillas} for the dual version of this result.
Clearly if the antipode of a Hopf algebra $B$  is bijective, then $\lgldim B = \rgldim B$.

We apply the estimate \eqref{eq:regular-gral} to a cleft exact sequence of Hopf algebras (all with bijective antipode)
$A\xhookrightarrow[]{\iota} C \xrightarrowdbl[]{\pi} B$ where we assume that $\dim A < \infty$.
First, we record

\begin{lemma}\label{lema:regular-ss} In the situation above, if $A$ is semisimple and $B$ is regular,
then $C$ is regular and $\gldim C \leq \gldim B$. \qed
\end{lemma}

Now  any finite-dimensional Hopf algebra is Frobenius, thus either $A$ is semisimple, or
$\gldim  A = \infty$.  The latter case is dealt with the following result, that relies on a theorem by Schneider.

\begin{lemma}\label{lemma:regular-gral} 
A Hopf algebra $K$ having a finite-dimensional normal non-semisimple Hopf subalgebra $L$ is not regular.
\end{lemma}

\pf  By \cite[Theorem 2.1]{sch}, $K$ is free as $L$-module. Therefore
\begin{align*}
\infty = \gldim L \overset{\heartsuit}{=}\pd \ku_L  \overset{\diamondsuit}{\leq}  \pd \ku_K + \pd K_L = \pd \ku_K,
\end{align*}
that is, $K$ is not regular. Here $\heartsuit$ is e.g. by \cite[2.4]{lorenz-lorenz} and $\diamondsuit$  by \cite[7.2.1]{mcconell-robson}.
\epf

Let  $\varGamma \leq G$.
Turning to the cleft extension $(\F^{\varGamma})$, we get the following characterization.

\begin{theorem}\label{thm:regular} Assume that  $\alg(\varGamma)$ is Noetherian. Then 
$\alg(\varGamma)$ is regular if and only if $H_{\varepsilon}$ is semisimple; in this case,
$\gldim\alg(\varGamma) \leq h$ where  $h$ is the Hirsch number of $\varGamma$. \qed
\end{theorem}

\pf By \cite{larson-radford}, $H_{\varepsilon}$ is  semisimple iff $H_{\varepsilon}^{*}$ is  so. 
If $H_{\varepsilon}^{*}$ is  not semisimple, then $\gldim \alg(\varGamma) = \infty$ 
by  Lemma \ref{lemma:regular-gral}. 
Otherwise Lemma \ref{lema:regular-ss} implies that $\gldim\alg(\varGamma) \leq \gldim \ku \varGamma$. When $\alg(\varGamma)$ is Noetherian,  
$\varGamma$ is polycyclic-by-finite by Theorem \ref{thm:Noetherian}, thus 
 $\gldim \ku \varGamma \leq h$ by \cite[7.5.6]{mcconell-robson}. 
\epf

\section{Examples}\label{section:examples}
\subsection{First considerations}
In this Section we discuss examples of Hopf algebras $\alg(\varGamma)$ as defined in Section \ref{sec:hopf-systems}, see \eqref{eq:def-HGamma}, and more specifically as assumed in Section \ref{sec:finite-by-coco}.
Once a substantial list of examples is obtained, to have a comprehensive picture one needs to address the following questions:

\begin{question}\label{question:quantization} Given an algebraic group  $G$  with algebra of functions $\Oc(G)$,
find all Hopf algebras $H$ with a central Hopf subalgebra $A \simeq \Oc(G)$
such that the extension
$A \overset{\iota}{\hookrightarrow} H \overset{p_{\varepsilon}}{\twoheadrightarrow} H_{\varepsilon}$
is cleft and  $\dim H_{\varepsilon}< \infty$. 
\end{question}

Question \ref{question:quantization} can be rephrased as follows:

\begin{question}\label{question:quantization2} 
Given an algebraic group  $G$  and a   Hopf algebra $\ug$, $\dim \ug < \infty$, 
find all cleft extensions
$\Oc(G) \overset{\iota}{\hookrightarrow} H {\twoheadrightarrow} \ug$
such that $\iota(\Oc(G))$ is central in $H$.
\end{question}

As is known \cite{ad,hofstetter,majid}, such an extension can be described by a collection $(\rightharpoonup, \sigma, \rho, \tau)$ made up of a weak action, a cocycle, a weak coaction and a cycle that satisfies a long set of axioms; 
 centrality  slightly simplifies the requirements (e.g. the weak action is trivial) but 
otherwise this situation seems to be difficult to handle.

\medbreak
Then we also need information on the following classical problem.

\begin{question}\label{question:fingen-virtuslly-nilpotent-subgroup}
 Given an algebraic group  $G$ (that admits a Noetherian Hopf algebra $H$ as in Question \ref{question:quantization}), 
 describe its subgroups, in particular those that are
finitely generated nilpotent-by-finite, or polycyclic-by-finite.
\end{question}

Towards this question, it is worth recalling the following celebrated result.

\begin{theorem}\label{thm:auslander-swan} \cite{auslander,swan}.
Any polycyclic-by-finite group is linear. \qed
\end{theorem}

\subsection{Quantum algebras of functions}\label{subsec:quantum-alg-functions}  
Let $G$ be a semisimple simply connected algebraic group  with Lie algebra $\g$ and algebra of functions $\Oc(G)$.
Let $\ell$ be an odd integer (prime to 3 if $G$ has a component of type $G_2$), 
and $\epsilon$ a primitive $\ell$-th root of unity. 
Recall the quantized algebra  of functions $\Oc_\epsilon(G)$, see e.g. \cite{deconcini-lyuba}, and the small quantum group
$\ug_\epsilon(\g)$.
There is an exact sequence of Hopf algebras 
\begin{align*}
\xymatrix{ \Oc(G)  \ar@{^{(}->}[r]  & \Oc_\epsilon(G) \ar@{->>}[r]^{\pi}  & {\ug_\epsilon(\g)}^{*}}
\end{align*}
which is cleft by \cite[3.4.3]{schauen-schn}. By  Theorem \ref{thm:auslander-swan}
any polycyclic-by-finite group $\varGamma$ is a subgroup of a suitable $G$ and so gives rise to a Hopf algebra
$\alg(\varGamma)$. 

More examples are given by the Hopf algebra
quotients of the   quantized algebras  of functions classified in \cite{ag-compo}.

\subsection{Pointed Hopf algebras}\label{subsec:deCKP-quantum}  
Recall that a Hopf algebra is \emph{pointed} if every simple comodule is one-dimensional. See \cite{as-cambr}
for this notion and \cite{a-leyva} for the related notion of a Nichols algebra. We start by  a result needed later.

\begin{theorem}\label{thm:masuoka} \cite[1.3]{masuoka} Let $\pi: U \to \ug$ be a surjective map of Hopf algebras. If $U$ is pointed,
then the $\ug$-comodule algebra $U$, with coaction  $(\id \otimes \pi)\Delta$,  is cleft. \qed
\end{theorem}

Let  $(V, c)$ be a braided vector space. Then the tensor algebra $T(V)$ is a braided graded Hopf algebra.
A \emph{pre-Nichols algebra} of $V$ is a  factor of $T(V)$ 
by a graded Hopf ideal  supported in degrees $\geqslant 2$. The Nichols algebra of $V$ is  
$\toba(V) = T(V) / \cJ(V)$, where $\cJ(V)$ is the maximal Hopf ideal among those. Thus any pre-Nichols algebra $\cB$ of $V$ 
lies between $T(V)$ and $\toba(V)$. 

\medbreak We refer to
\cite{aa-diag-survey,aay, angiono-spn} for details on the following material.
Let  $(V, c)$ be a  braided vector space of diagonal type with braiding matrix 
$\bq = (q_{ij} ) \in \big( \ku^\times \big)^{\I \times \I}$, where $\I = \{1, \dots, \theta\}$, $\theta \in \N$. 
Assume that the Nichols algebra 
$\toba_{\bq} \coloneqq\toba(V)$ has finite dimension, thus $\bq$ belongs to the classification in \cite{H-classif RS}.
As shown in \cite{angiono-spn}, $(V,c)$, i.e., the matrix $\bq$, gives rise to the following data:

\begin{itemize}[leftmargin=*]\renewcommand{\labelitemi}{$\circ$}

\medbreak
\item The distinguished pre-Nichols algebra  $\dpn_{\bq}$. 

\medbreak
\item The Hopf algebra $U_\bq$ 
(the  Drinfeld double of the bosonization of $\dpn_{\bq}$). 
 
 \medbreak
 \item The  subalgebra $\mathcal Z_{\bq}$ of $U_{\bq}$ 
 as modified in \cite[\S 4.5]{aay}.
\end{itemize}

The class of Examples of this Subsection arises from the following result.
\begin{theorem} 
If $\bq$ satisfies the technical condition \cite[(4.26)]{aay}, then  $\mathcal Z_{\bq}$ is a central Hopf subalgebra of $U_{\bq}$,
$M_{\bq}  \coloneqq \Alg(\mathcal Z_{\bq}, \ku)$ is a solvable algebraic group,
$\widetilde{\ug}_{\bq} = U_{\bq}/ U_{\bq} \mathcal Z_{\bq}^+$ is \fd, 
and the  exact sequence   of Hopf algebras 
$\xymatrix{\mathcal Z_{\bq} \ar@{^{(}->}[r]  & U_{\bq} \ar@{->>}[r]  & \widetilde{\ug}_{\bq}}$
is cleft.
\end{theorem}

\pf  This follows from \cite[Theorem 33 \& Remark 11]{angiono-spn}, see also the discussion in \cite[\S 4.5]{aay}.
Theorem \ref{thm:masuoka} implies the cleftness of the exact sequence.
\epf

Let $\g$ be a semisimple Lie algebra and $q$ a root of 1 of odd order, coprime with $3$ if  $\g$ has an ideal of type $G_2$.
Then there is a suitable $\bq$ such that  $U_{\bq}$ is isomorphic to
the De Concini-Kac-Procesi quantized enveloping algebra $U_q(\g)$, 
see  \cite{deconcini-kac,deconcini-kac-procesi,deconcini-procesi}.
But this class of examples covers also quantum supergroups and more, see \cite{aa-diag-survey,aay}.

\subsection{Bosonizations}\label{subsec:examples-bosonization}
See e.g. \cite{as-cambr}
for this notion due to Radford and interpreted categorically by Majid.
The category of Yetter-Drinfeld modules over a Hopf algebra $K$ is denoted by $\yd{K}$.

\medbreak
We fix  a finite group $L$, $V \in \yd{\ku L}$ such that $\dim \toba(V) < \infty$ and set $K = \toba(V) \# \ku L$.
We consider three different settings for pairs $(H, A)$  with $H$ Noetherian, $A$ central in $H$, the extension
$(\E)$ cleft and $H_{\varepsilon} \simeq K$.

\medbreak
\subsubsection{Blowing-up the group} \label{subsubsec:blow-group}
We assume that there exist a group $\Lambda$, a  surjective map of groups $\wp: \Lambda \to L$
and a map $\varsigma:\supp V \to \Lambda$ such  that  $V \in \yd{\ku\Lambda}$ with 
the grading induced by $\varsigma$ and the action induced by $\wp$. Set  $\Lambda_0 \coloneqq \ker \wp$.
We claim that the following  sequence is exact and cleft:
\begin{align*}
\xymatrix{ A' \coloneqq \ku \Lambda_0 \  \ar@{^{(}->}[rr] & & 
H' \coloneqq \toba(V) \# \ku \Lambda \ar@{->>}[rr]^{\qquad \id \# \wp} & &K.}
\end{align*} 
Indeed tensoring the exact sequence of vector spaces $(\ku \Lambda) \ku \Lambda_0^+ \hookrightarrow \ku \Lambda \twoheadrightarrow \ku L$ with $\toba(V)$, we conclude that $\ker (\id \# \wp) = H'(A')^+$.  It is also easy to see
that $A' \overset{\diamond}{\subseteqq} (H') ^{\operatorname{co} \id \# \wp}$. Now $H'$ is $K$-cleft 
by Theorem \ref{thm:masuoka}. Hence we have equality in $\diamond$ by a dimension argument. 
Finally we need to check  whether $\Lambda_0$  is central in $H'$ in each example.

\subsubsection{Blowing-up the algebra}\label{subsubsec:blow-algebra} 
Here we suppose that there  exists a pre-Nichols algebra $\cB$  of $V$ with projection  $\pi:\cB \to \toba(V)$
such that $\cB^{\operatorname{co} \pi}$ is a (usual) Hopf subalgebra and  the   following  sequence is exact and cleft:
\begin{align*}
\xymatrix{ A^{\prime\prime} \coloneqq \cB^{\operatorname{co} \pi} \  \ar@{^{(}->}[rr] & & 
H^{\prime\prime}  \coloneqq \cB \# \ku L \ar@{->>}[rr]^{\pi \# \id} & &K.}
\end{align*}
We also assume that $A^{\prime\prime}$ is central in $H^{\prime\prime}$.

\subsubsection{Blowing-up both} \label{subsubsec:blow-both}
We assume that there  exist a group $\Lambda$ and a pre-Nichols algebra $\cB$  as 
previously such that  $\cB^{\operatorname{co} \pi}\# \ku \Lambda_0$ is a central Hopf subalgebra 
and  the   following  sequence is exact and cleft:
\begin{align*}
\xymatrix{ A \coloneqq \cB^{\operatorname{co} \pi} \# \ku \Lambda_0\  \ar@{^{(}->}[rr] & & 
H \coloneqq \cB \# \ku \Lambda \ar@{->>}[rr]^{\pi \# \wp} & &K.}
\end{align*}

\subsubsection{Example: the enveloping group of a rack}\label{subsubsec:FK-enveloping group}
For simplicity we fix $n \in \I_{3,5}$ and assume that $L = \sn$ and that $V \in \yd{\ku \sn}$ 
has $\toba(V)  \simeq \FK_{n}$ (the Fomin-Kirillov algebra of rank $n$),  see e.g. \cite{a-leyva}. 
Then $V$ has support $X = \Oc^n_2$ (the conjugacy class of transpositions in $\sn$);
this is a rack with operation $\trid$ given by conjugation.
The  enveloping group of $X$ is
\begin{align*}
\Gx &:= \langle e_x: x\in X\vert e_x\, e_y = e_{x\trid y}\, e_x, x,y \in X\rangle.
\end{align*}
Then there exist \begin{enumerate*}[label=\rm{(\roman*)}, itemjoin={,\quad}] 
\item a unique morphism of groups $\wp: \Gx \to \sn$ such that
$\wp(e_x) = x$, $x\in X$ \item a realization of  $V$ in $\yd{\Gx}$,
\end{enumerate*} inducing a Hopf algebra map
$\xymatrix{ H' \coloneqq \toba(V) \# \ku \Gx \ar@{->>}[r]^{\id \# \wp}  &\toba(V) \# \ku \sn.}$
It can be shown that $z \coloneqq e_x^2 =e_y^2$, for any $x,y \in X$ and that 
$Z \coloneqq\ker \wp = \langle z\rangle \simeq \Z$. 
Now it is clear that  $e_x^2 \cdot v = v$ for every $v \in V$, hence $Z$ is central in  
$H'$. Thus we have a cleft exact  sequence
\[
\xymatrix{ A' \coloneqq \ku Z   \ar@{^{(}->}[rr]  & &
	H'  \ar@{->>}[rr]^{\id \# \wp \qquad} && K = \toba(V) \# \ku \sn.}\]
Now it is well-known that
\begin{align*}
\Alg (\ku Z, \ku) &\simeq \widehat{Z} \simeq \kut, &
K^* &\simeq \toba(W) \# \ku^{\sn},
\end{align*}
where $W \simeq V^*$ linearly. Hence for any $\varGamma \leq \kut$, we have an extension
\begin{align*}
\xymatrix{ \toba(W) \# \ku^{\sn}  \ar@{^{(}->}[r]  & 
\alg(\varGamma)  \ar@{->>}[r] & \ku \varGamma.}
\end{align*}

\subsubsection{Example: (finite) quantum linear spaces}\label{subsubsec:FK-QLS}
The input for this example is a matrix $\bq = (q_{ij})_{i, j \in \I_{\theta}}$ whose entries are roots of 1  that satisfy
\begin{align*}
q_{ij}q_{ji} &= 1, & i\neq j &\in \I_{\theta}, & N_i &\coloneqq \ord q_{ii} >1, & i  &\in \I_{\theta}.
\end{align*}
Let $M_i$ be the least common multiple of $\{\ord q_{ij}: j \in \I_{\theta} \}$, for every $i\in \I_{\theta}$.
Thus $N_i$ divides $M_i$.
To this input we attach:
\begin{itemize}[leftmargin=*]\renewcommand{\labelitemi}{$\diamond$}
\item The finite abelian group $L = \langle g_{1}\rangle \oplus \cdots \oplus \langle g_{\theta} \rangle$ where 
$\ord g_{i} = M_i$, $i\in \I_{\theta}$.

\medbreak
\item The characters $\chi_{1}, \dots \chi_{\theta} \in \widehat L$, given by $\chi_{j}(g_{i}) = q_{ij}$, $i, j \in \I_{\theta}$.

\medbreak
\item The free abelian group $\Lambda$ of rank $\theta$, with basis $\gb_{1}, \dots \gb_{\theta}$. 
The kernel of the map $\wp: \Lambda \to L$ given by $\wp(\gb_i) = g_i$, $i\in \I_{\theta}$, is denoted by $\Lambda_0$; that is,
$\Lambda_0 = \langle \gb_{1}^{M_1}, \dots \gb_{\theta}^{M_{\theta}}\rangle$.

\medbreak
\item The characters $\chib_{1}, \dots \chib_{\theta} \in \widehat \Lambda$, given by $\chib_{j}(\gb_{i}) = q_{ij}$, $i, j \in \I_{\theta}$.

\medbreak 
\item
A vector space $V$ with a basis $x_{1}, \dots x_{\theta}$, realized in $\yd{\ku L}$ and $\yd{\ku \Lambda}$
by
\begin{align*}
\delta(x_i)  &= g_i \otimes x_i, & g\cdot x_i &= \chi_{i}(g) x_i, &
{\boldsymbol \delta}(x_i)  &= \gb_i \otimes x_i, & \gb\cdot x_i &= \chib_{i}(x) x_i, 
\end{align*}
for  $i\in \I_{\theta}$, $g \in L$, $\gb \in \Lambda$. Here $\delta: V \otimes \ku L$ 
and ${\boldsymbol\delta}: V \otimes \ku \Lambda$ are the coactions.

\medbreak
\item The algebra $\cB = \ku \langle x_{1}, \dots ,x_{\theta} \vert x_ix_j - q_{ij} \, x_j x_i, i\neq j \in \I_{\theta}\rangle$.
\end{itemize} 
Then the following facts hold:

\begin{itemize}[leftmargin=*]
\item  $\toba(V)$ and $\cB$ are  both Hopf algebras in  $\yd{\ku L}$ and $\yd{\ku \Lambda}$.
The natural projection $\pi$ from $\cB$ to the Nichols algebra $\toba(V)$ induces an isomorphism 
$\toba(V) \simeq\cB / \langle x_{1}^{N_1}, \dots, x_{\theta}^{N_\theta}\rangle$.

\medbreak
\item  The subalgebra  of $\cB$ generated by $x_{1}^{N_1}, \dots, x_{\theta}^{N_\theta}$ coincides with
$\cB^{\operatorname{co} \pi}$  and we have an exact sequence of Hopf algebras
\begin{align}\label{eq:exact-sequence-QLS}
\xymatrix{ A \coloneqq \cB^{\operatorname{co} \pi} \# \ku \Lambda_0\  \ar@{^{(}->}[rr] & & 
H \coloneqq \cB \# \ku \Lambda \ar@{->>}[rr]^{\pi \# \wp} & &K.}
\end{align}
\end{itemize} 

\begin{lemma}
If $N_i = M_i$, for all $i\in \I_{\theta}$, then  $A$ is central in $H$.
\end{lemma}

\pf The algebra $H$ is generated by $x_{1}, \dots ,x_{\theta}$, $\gb_{1}^{\pm 1}, \dots ,\gb_{\theta}^{\pm 1}$ with relations
\begin{align*}
 x_ix_j &= q_{ij} \, x_j x_i, & i\neq j &\in \I_{\theta}, &
  \gb_ix_j &= q_{ij} \, x_j \gb_i, & 
   \gb_i \gb_j &=  \gb_ j \gb_i, & i, j &\in \I_{\theta}.
\end{align*}
The subalgebra $A$ is generated by $x_{1}^{N_1}, \dots, x_{\theta}^{N_\theta}$, 
$\gb_{1}^{\pm N_1}, \dots ,\gb_{\theta}^{\pm N_\theta}$. Now
\begin{align*}
x_i^{N_i}x_j &= q_{ij}^{N_i} \, x_j x_i^{N_i}, &
\gb_i^{N_i} x_j &= q_{ij}^{N_i} \, x_j \gb_i^{N_i}, & 
 \gb_ix_j^{N_j} &=\begin{cases}
q_{ii}^{N_i} \, x_i^{N_i} \gb_i, &\text{ if } i=j;\\
q_{ji}^{-N_j} \, x_j^{N_j} \gb_i &\text{ if } i \neq j,
  \end{cases} 
\end{align*}
for all $i,j \in \I_{\theta}$.
This implies the Lemma. \epf

Now  the comultiplication of $A$ is determined by 
\begin{align*}
\Delta (\gb_i^{N_i}) &= \gb_i^{N_i} \otimes \gb_i^{N_i}, &
\Delta (x_i^{N_i}) &= x_i^{N_i} \otimes 1 + \gb_i^{N_i} \otimes x_i^{N_i}, & i \in \I_{\theta}.
\end{align*}
Therefore
$G \coloneqq \Alg (A, \ku)$ is isomorphic to $\mathbb{B}^{\theta}$, where $\mathbb{B}$ is the Borel subgroup of $SL(2, \ku)$. 
In conclusion,  for any $\varGamma \leq \mathbb{B}^{\theta}$, we have an extension
\begin{align*}
\xymatrix{ K^{*}  \ar@{^{(}->}[r]  & 
\alg(\varGamma)  \ar@{->>}[r] & \ku \varGamma.}
\end{align*}

\subsection{Twisting $\alg(\varGamma)$}\label{subsec:twisting-datum}
Let  $A$ be  a central Hopf subalgebra of a Noetherian Hopf algebra $H$,
so that one has the exact sequence  
$A \overset{\iota}{\hookrightarrow} H \overset{p_{\varepsilon}}{\twoheadrightarrow} H_{\varepsilon}$
of Hopf algebras.
 Given  $\varGamma \leq G = \Alg(A, \ku)$, let  $\alg(\varGamma)$ be as in \eqref{eq:def-HGamma}.
Let 
$\sigma:H_{\varepsilon} \ot H_{\varepsilon}\to \ku$ be a cocycle, cf. \S \ref{subsec:cocycles-twists}. Then 
$\widetilde{\sigma} \coloneqq \sigma(p_{\varepsilon} \otimes p_{\varepsilon}): H \otimes H \to \ku$
is a cocycle.  Then
\begin{align*}
	z \cdot_\sigma h &= zh, &h \cdot_\sigma z &= hz, & z \in A, \, h &\in H. 
\end{align*}
Hence $A$ is central in $H_{\widetilde{\sigma}}$, 
Assume that $H$ is a finitely generated $A$-module. Then
$H_{\widetilde{\sigma}}$, being a finitely generated $A$-module, is Noetherian 
and we have an exact sequence of Hopf algebras
$\xymatrix{ A  \ar@{^{(}->}[r]^{\iota} & H_{\widetilde{\sigma}} \ar@{->>}[r]^{p_{\varepsilon}} & (H_{\varepsilon})_{\sigma}}$.
Thus we may repeat the constructions 
and get a Hopf algebra that we denote $\alg_{\sigma}(\varGamma)$.

\medbreak
The twist $F = {}^t\sigma$ for $H_{\varepsilon}^{*} \simeq C(\varepsilon)$ is naturally
a twist for $\alg(\varGamma)$ and for $H^\circ$. See Example \ref{exa:twist-graded} for twists of graded algebras.

\begin{prop}\label{prop:twist} $\alg_{\sigma}(\varGamma) \simeq \alg(\varGamma)^F$.
\end{prop}

\noindent\emph{Proof.} Since the subcoalgebras $C(\kappa)$ of $\alg(\varGamma)$ are $C(\varepsilon)$-bimodules, then $C(\kappa)^F$ is a subcoalgebra of $\alg(\varGamma)^F$, for all $\kappa \in \Gamma$.
Notice that  $(H_{\widetilde{\sigma}})._{\widetilde{\sigma}}\mathfrak{M}_\kappa = H\mathfrak{M}_\kappa$ for all $\kappa \in \Gamma$. Hence the subcoalgebra $C_{\widetilde{\sigma}}(\kappa) :=  (H_{\widetilde{\sigma}}/(H_{\widetilde{\sigma}})._{\widetilde{\sigma}}\mathfrak{M}_\kappa)^*$ of $(H_{\widetilde{\sigma}})^\circ$ is contained in $(H^\circ)^F$ and it coincides with the subcoalgebra $C(\kappa)^F$. Therefore
\begin{align*}
\alg(\varGamma)^F = \oplus_{\kappa \in \varGamma} C(\kappa)^F = \oplus_{\kappa \in \varGamma}  
C_{\tilde{\sigma}}(\kappa) = \alg_{\sigma}(\varGamma). \qed
\end{align*}

Proposition \ref{prop:twist} has the following application. Let $\ug$ be a finite-dimensional pointed Hopf algebra.
In all known examples, the graded Hopf algebra $\gr \ug$ with respect to the coradical filtration is of the form
$\toba(V) \# \ku L$ as in \S \ref{subsec:examples-bosonization} and $\ug \simeq (\gr \ug)_{\sigma}$ for a suitable cocycle $\sigma$.  Then for any pair $(H, A)$  with $H$ Noetherian, $A$ central in $H$, the extension
 $(\E)$ cleft and $H_{\varepsilon} \simeq \gr \ug$, the pair $(H_{\widetilde{\sigma}}, A)$   has analogous properties and
 $(H_{\varepsilon})_{\sigma} \simeq \ug$ .

\subsection*{Acknowledgements} 
N.A. thanks Serge Skryabin for a helpful email exchange and Juan Cuadra for interesting virtual conversations in the first semester of 2020 that were the genesis of this paper.


\begin{thebibliography}{11}

\bibitem{a-can} N. Andruskiewitsch.
{\it  Notes on extensions of Hopf algebras},
Canad. J. Math.  \textbf{48}, 3--42 (1996).

\bibitem{a-leyva} \bysame. \emph{An Introduction to Nichols Algebras}. In Quantization, Geometry and Noncommutative Structures in Mathematics and Physics. 
A. Cardona et al., eds., pp. 135--195, Springer-Nature (2017).


\bibitem{aa-diag-survey} N. Andruskiewitsch and I. Angiono. \emph{On finite dimensional Nichols algebras of diagonal type}. Bull. Math. Sci. \textbf{7}, 353--573 (2017). 

\bibitem{aay} N. Andruskiewitsch, I. Angiono, M. Yakimov. \emph{Poisson orders on large quantum groups}.
\texttt{arXiv:2008.11025} (2020).

\bibitem{AC} N. Andruskiewitsch, J. Cuadra. \emph{On the structure of (co-Frobenius) Hopf algebras}. J.
Noncommut. Geom. \textbf{7}  83--104 (2013).

\bibitem{ace} N. Andruskiewitsch, J. Cuadra, P. Etingof.
\emph{On two finiteness conditions for Hopf algebras with nonzero integral}.
Ann. Sc. Norm. Super. Pisa  \textbf{14}, 401--440 (2015).


\bibitem{andrusdasca} N. Andruskiewitsch and S. D\v{a}sc\v{a}lescu.  
\emph{Co-Frobenius Hopf algebras and the coradical filtration}.
 Math. Z. \textbf{243},  145--154 (2003).


\bibitem{ad} N. Andruskiewitsch and J. Devoto. \emph{Extensions of Hopf algebras}.
Algebra i Analiz, \textbf{7} (1)  22–61 (1995); St. Petersburg Math. J. \textbf{7} (1)  17--52 (1995).





\bibitem{ag-compo}  N. Andruskiewitsch and  G. A. Garc\'\i a. 
\emph{Finite subgroups of a simple quantum group}. Compositio Math. \textbf{145}, pp. 476--500 (2009).



\bibitem{as-cambr} N. Andruskiewitsch, H.-J. Schneider. 
\emph{Pointed Hopf algebras}. In Recent developments in Hopf algebras Theory, 
MSRI Publ. \textbf{43} 1--68, Cambridge Univ. Pr. (2002).

\bibitem{angiono-spn}  I. Angiono. \emph{Distinguished pre-Nichols algebras}.
Transform. Groups {\bf{21}}, 1--33 (2016). 


\bibitem{auslander} L. Auslander. \emph{On a problem of Philip Hall}. 
Ann. Math. (2) \textbf{86}, 112--116 (1967).


\bibitem{bcm} R. J. Blattner, M. Cohen and S. Montgomery.
\emph{Crossed products and inner actions of Hopf algebras}. Trans.
Amer. Math. Soc.    \textbf{298},  671--711 (1986).



\bibitem{bc}
K. A. Brown and M. Couto. \emph{Affine commutative-by-finite Hopf algebras}.  J. Algebra \textbf{573}:1, 56--94 (2021).

\bibitem{bcj}
K. A. Brown, M. Couto, and A. Jahn. \emph{The finite dual of commutative-by-finite Hopf algebras}.
Glasgow Math. J. \textbf{65}, 62--89. \texttt{doi:10.1017/S0017089522000052}.


\bibitem{byott} N. P. Byott. \emph{Cleft extensions of Hopf algebras}. J. Algebra 1\textbf{57}, 405--429 (1993).
\newline  \emph{Cleft extensions of Hopf algebras II}. Proc. London Math. Soc. (3)  \textbf{67}, 277--304 (1993).



\bibitem{chirvasitu-kasprzak} A. Chirvasitu and P. Kasprzak.
\emph{On the Hopf (co)center of a Hopf algebra}.
J. Algebra \textbf{464}, pp. 141--174 (2016).


\bibitem{dasca-nasta-torrecillas} 
S. D\u{a}sc\u{a}lescu, C.  N\u{a}st\u{a}sescu, B. Torrecillas. 
\emph{Homological dimension of coalgebras and crossed coproducts}. K-Theory \textbf{23}, pp. 53--65 (2001). 

\bibitem{deconcini-kac} C. De Concini and V. Kac. \emph{Representations of quantum groups at roots of $1$}.
 Progr. Math. \textbf{92},  471--506, Birkh\"auser, Boston (1990).

\bibitem{deconcini-kac-procesi} C. De Concini, V. Kac, C. Procesi.
\emph{Quantum coadjoint action}. J.~Amer. Math. Soc. \textbf{5},  151--189 (1992).


\bibitem{deconcini-lyuba} C. De Concini and  V. Lyubashenko. \emph{Quantum function algebra at roots of 1}.
Adv. Math. \textbf{108}, No. 2, 205--262 (1994).


\bibitem{deconcini-procesi} C. De Concini and C. Procesi. 
\emph{Quantum groups}. In: D-modules, representation theory and quantum groups, 
Lecture Notes in Math. \textbf{1565} 31--140, Springer (1993).


\bibitem{doi1} Y. Doi. \emph{On the structure of relative Hopf modules}.  Comm. Algebra \textbf{11}, 243--255 (1983).


\bibitem{doi} \bysame. \emph{Braided bialgebras and quadratic bialgebras}, Comm. Algebra {\bf 21}, 1731--1749  (1993).

\bibitem{doi-takeuchi}  Y. Doi and M. Takeuchi.  \emph{Cleft comodule algebras for a bialgebra}. 
Comm. Algebra    \textbf{14}      801--818 (1986).

\bibitem{EGNO}
P.~Etingof, S.~Gelaki, D.~Nikshych, V.~Ostrik, {\em  Tensor categories}.
Mathematical Surveys and Monographs,
\textbf{205}, American Mathematical Society (2015).

\bibitem{formanek-jategaonkar} E. Formanek and A. V. Jategaonkar.
\emph{Subrings of Noetherian rings}.
Proc. Amer. Math. Soc. \textbf{46}, 181--186 (1974).

\bibitem{gromov} M. Gromov.
\emph{Groups of polynomial growth and expanding maps. Appendix by Jacques Tits}. 
Publ. Math., Inst. Hautes Étud. Sci. \textbf{53}, 53--78 (1981).

\bibitem{hall} P. Hall. \emph{Finiteness conditions for soluble groups}. Proc. London Math. Soc. (3) \textbf{4},
419--436 (1954). 

\bibitem{H-classif RS} I. Heckenberger.  \emph{Classification of arithmetic root systems}. Adv. Math. \textbf{220}, 59--124 (2009).

\bibitem{hofstetter} I. Hofstetter. \emph{Extensions of Hopf algebras and their cohomological description}.
J. Algebra \textbf{164} 264--298  (1994).


\bibitem{ivanov}  S.V. Ivanov. \emph{Group rings of Noetherian groups}. 
Math. Notes \textbf{46},  929--933 (1989); translation from Mat. Zametki \textbf{46}, 61--66 (1989).



\bibitem{kac-extension} Kac, G.,  \emph{Extensions of groups to ring groups}. Math. USSR. Sb. \textbf{5}, 451--474 (1968).

\bibitem{krause-lenagan} G. Krause and T. Lenagan.  \emph{Growth of algebras and Gelfand-Kirillov dimension}. Revised ed.
Graduate Studies in Math. \textbf{22}. Amer. Math. Soc. (2000).


\bibitem{li-liu} K. Li, G. Liu. 
\emph{The Finite Duals of Affine Prime Regular Hopf Algebras of GK-Dimension One}.
AIMS Mathematics,  \textbf{8} (3), 6829--6879 (2023).

\bibitem{larson-radford} R. G. Larson, D. Radford. \emph{Finite dimensional cosemisimple Hopf algebras in characteristic 0 are semisimple}. J. Algebra \textbf{117}, 267--289 (1988).

\bibitem{lorenz-lorenz} M. Lorenz, M. Lorenz. 
\emph{On crossed products of Hopf algebras}. Proc. Amer. Math. Soc. \textbf{123}, 33--38 (1995).


\bibitem{majid} S. Majid. \emph{Crossed products by braided groups and bosonization}. J. Algebra \textbf{163}, 165--190  (1994).

\bibitem{masuoka} A. Masuoka. \emph{On Hopf algebras with cocommutative coradicals}. 
J. Algebra \textbf{144}, 415--466  (1991).

\bibitem{mcconell-robson} J.C. McConnell and J.C. Robson. \emph{Noncommutative Noetherian rings. With the cooperation of L. W. Small}.  Reprinted with corrections from the 1987 original. Providence, RI: Amer. Math. Soc. (AMS) (2001).

\bibitem{molnar} R. K. Molnar. \emph{A commutative Noetherian Hopf algebra over a field is finitely generated}.
Proc. Amer. Math. Soc. \textbf{51}, 501--502 (1975).




\bibitem{montgomery} S. Montgomery.
\emph{Hopf Algebras and their Actions on Rings}. AMS (1993), CMBS {\bf 82}.

\bibitem{nastasescu-fvo}
C. Nastasescu and F. Van Oysteyen. \emph{On strongly graded rings and crossed products}. 
Commun.  Algebra \textbf{10}, 2085--2106 (1982).

\bibitem{nz}  W.D. Nichols and M. B. Z\"oller. \emph{A Hopf
algebra freeness Theorem}. Amer. J.  Math. \textbf{111},  381--385 (1989).  

\bibitem{rad-libro} Radford, D. E. \emph{Hopf algebras}. Series on Knots and Everything 49. 
Hackensack, NJ: World Scientific. xxii, 559 p.  (2012).


\bibitem{schauen-schn} P. Schauenburg, H.-J. Schneider. \emph{Galois type extensions and Hopf algebras}. (2004).


\bibitem{sch} H.-J. Schneider. \emph{Some remarks on exact sequences of quantum groups}. Comm. Algebra   \textbf{21} (9) 3337--3358 (1993).

\bibitem{sch2}  \bysame.   \emph{Zerlegbare Erweiterungen affiner Gruppen}. 
J. Algebra  \textbf{66}, 569--593 (1980) .

\bibitem{sch3} \bysame.   \emph{Zerlegbare Untergruppen affiner Gruppen}.   Math. Ann.  \textbf{255},    139--158 (1981).



\bibitem{singer} W. Singer.  \emph{Extension theory for connected Hopf algebras}.  J. Algebra \textbf{21}, 1--16 (1972).

\bibitem{swan} R. G. Swan. \emph{Representations of polycyclic groups}. 
Proc. Amer. Math. Soc. \textbf{18}, 573--574 (1967).



\bibitem{takeuchi} M. Takeuchi. \emph{Matched pair of groups and bismash product of Hopf algebras}.
Comm. Algebra \textbf{9}, 841--882 (1981).

\bibitem{takeuchi-quotient} \bysame. \emph{Quotient spaces for Hopf algebras}. Comm. Algebra  \textbf{22},  2503--2523 (1994).

\bibitem{tits} J. Tits. \emph{Free subgroups in linear groups}.
J. Algebra \textbf{20}, 250--270 (1972).



\bibitem{ulbrich} K. H. Ulbrich. \emph{Vollgraduierte Algebren}. Abh. Math. Semin. Univ. Hamb. \textbf{51}, 136--148 (1981).

\end{thebibliography}
\end{document}